\documentclass[10pt]{article}

\usepackage[latin1]{inputenc}
\usepackage{url, amsmath,enumerate,fancyhdr,amssymb, amsthm, 
pstricks, epsf, lscape, ifpdf, graphicx, float}
\usepackage{algorithm}
\usepackage{color}
\usepackage{algpseudocode} 
\algnewcommand\algorithmicinput{\textbf{Input:}}
\algnewcommand\INPUT{\item[\algorithmicinput]}
\algnewcommand\algorithmicinitialization{\textbf{Initialization:}}
\algnewcommand\INITIALIZATION{\item[\algorithmicinitialization]}

\usepackage[margin=3cm]{geometry}

\mathchardef\mhyphen="2D

\newtheorem{Definition}{Definition}

\newtheorem{Example}{Example}
\newtheorem{Proposition}{Proposition}
\newtheorem{Lemma}{Lemma}
\newtheorem{Theorem}{Theorem}
\newtheorem{Corollary}{Corollary}
\newtheorem{Remark}{Remark}[section]
\newtheorem{Assumption}{Assumption}

\newcommand{\A}{\mathcal A}
\newcommand{\B}{\mathcal B}
\newcommand{\lin}{\operatorname{lin}}

\newcommand{\cl}{\operatorname{cl}}
\newcommand{\ri}{\operatorname{ri}}
\newcommand{\dir}{\operatorname{dir}}
\newcommand{\ldir}{\operatorname{ldir}}

\newcommand{\scr}{\scriptstyle}

\newcommand{\eps}{\epsilon} 
\newcommand{\bpx}{\begin{pmatrix}}
\newcommand{\epx}{\end{pmatrix}}
\newcommand{\bbx}{\begin{bmatrix}}
\newcommand{\ebx}{\end{bmatrix}}

\newcommand{\bdef}{\begin{Definition}} 
\newcommand{\commentout}[1]{}
\newcommand{\co}[1]{}

\newcommand{\coab}[1]{}

\newcommand{\nin}{\noindent}
\newcommand{\ti}{\times}
 
\newcommand{\pf}[1]{\vspace{.35cm} \nin {\bf Proof {#1} }}

\newcommand{\sym}[1]{{\cal S}^{#1}}
\newcommand{\psd}[1]{{\cal S}_+^{#1}}
\newcommand{\rad}[1]{\mathbb{R}^{#1}}

\newcommand{\eref}[1]{(\ref{#1})}

\newcommand{\R}{ {\cal R} }
\newcommand{\N}{ {\cal N} }
\newcommand{\M}{ {\cal M} }

\newcommand{\la}{\langle}
\newcommand{\ra}{\rangle}

\newcommand{\myint}{\operatorname{int}}

\newcommand{\beq}{\begin{equation}}
\newcommand{\eeq}{\end{equation}}
\newcommand{\beqa}{\begin{eqnarray}}
\newcommand{\eeqa}{\end{eqnarray}}
\newcommand{\ba}{\begin{array}}
\newcommand{\ena}{\end{array}}
\newcommand{\bac}{\begin{array}{ccccccccccc}}
\newcommand{\eac}{\end{array}}
\newcommand{\bprop}{\begin{Proposition}}
\newcommand{\eprop}{\end{Proposition}}

\newcommand{\np}{{\cal NP}}
\newcommand{\conp}{{\rm co \mhyphen}{\cal NP}}

\newcommand{\psdsysp}{(P_{SD}')}
\newcommand{\psdsysb}{(P_{\text{SD,bad}})}
\newcommand{\psdsysg}{(P_{\text{SD,good}})}

\newcommand{\beqast}{\begin{eqnarray*}}
\newcommand{\eeqast}{\end{eqnarray*}}
\newcommand{\benum}{\begin{enumerate}}
\newcommand{\eenum}{\end{enumerate}}
\newcommand{\bit}{\begin{itemize}}
\newcommand{\eit}{\end{itemize}}
\newcommand{\bth}{\begin{Theorem}}
\newcommand{\enth}{\end{Theorem}}
\newcommand{\ble}{\begin{Lemma}}
\newcommand{\ele}{\end{Lemma}}
\newcommand{\bex}{\begin{Example}}
\newcommand{\eex}{\end{Example}}
\newcommand{\bcor}{\begin{Corollary}}
\newcommand{\ecor}{\end{Corollary}}
\newcommand{\brem}{\begin{Remark}}
\newcommand{\erem}{\end{Remark}}
\newcommand{\bass}{\begin{Assumption}}
\newcommand{\eass}{\end{Assumption}}

\newcommand{\LRA}{\Leftrightarrow} 
\newcommand{\RA}{\Rightarrow} 
\newcommand{\LA}{\Leftarrow}

\newcommand\sdpspace{\hspace{.9cm}}

\newcommand{\val}{\operatorname{val}}
\newcommand{\aval}{\operatorname{aval}}

\newcommand{\newpc}{$\mathit{(P_c)\,}$}
\newcommand{\newdc}{$\mathit{(D_c)\,}$}

\newcommand{\pc}{\mathit{P_c}}
\newcommand{\dc}{\mathit{D_c}}

\newcommand{\sdpc}{\mathit{SDP_c}}
\newcommand{\sddc}{\mathit{SDD_c}}

\newcommand\emp{\emptyset}

\newcommand{\bsmx}{\begin{small} \begin{pmatrix}}
\newcommand{\esmx}{\end{pmatrix} \end{small}}

\title{\Large Bad semidefinite programs: they all look the same}
\author{G\'{a}bor  Pataki\thanks{Department of Statistics and Operations Research, University of North Carolina at Chapel Hill} \hspace{1cm} 
}

\begin{document}

\maketitle 

\begin{abstract} 

Conic linear programs, among them semidefinite programs,  often behave pathologically: 
the optimal values of the primal and dual programs may differ,
and may not be attained. We present a novel analysis of these pathological behaviors.
We call a conic linear system $\A x \leq_K b$ {\em badly behaved} if  the value of
$\sup \, \{ \la c, x \ra | \A x \leq_K b \}$ is finite but the  dual program has no solution with the same value for {\em some}   $c.$
We describe simple and intuitive geometric characterizations of 
badly behaved conic linear systems.
Our main motivation is the striking similarity of 
badly behaved semidefinite systems in the literature; 
we characterize such systems by certain {\em excluded matrices,} 
which are easy to spot in all published examples. 

We show how to transform semidefinite systems into a canonical form, which 
allows us to easily verify whether they are badly behaved. We prove several other structural results about badly behaved semidefinite systems; for example, we 
show that they are in $\np \cap \conp$ in the 
real number model of computing.
As a byproduct, we prove that all linear maps that act on symmetric matrices 
can be brought into a canonical form; this canonical form allows us
to easily check whether the image of the semidefinite cone under the given linear map is closed.

\co{and that the bad behavior of the transformed system is trivial to verify by inspection.}

\end{abstract}



{\em OR/MS subject classification:} Primary: convexity; secondary: programming-nonlinear-theory

\section{Introduction} \label{sect-intro}

Many  problems in engineering, combinatorial optimization, machine learning,  and related fields 
can be formulated as the primal-dual pair of conic linear  programs 
\begin{center}                                                                                              
$$                                                                                                          
\begin{array}{lrlcrlr}                                                                                      
         &   \sup  & \la c, x  \ra      &    \sdpspace & \inf   &   \la b, y  \ra &   \\ 
(\pc)  &   s.t.  & \A x \leq_K b        &   \sdpspace & s.t.  &   y \geq_{K^*} 0   & (\dc)     \\ 
         &         &                      &    \sdpspace &       &   \A^*y = c,  &    
\end{array}                                                                                              
$$   
\end{center}                                                                                                      
where 
$\A: X \rightarrow Y$ is a linear map between finite dimensional Euclidean spaces $X$ and $Y, \,$ 
$\A^*$ is its adjoint, 
$K \subseteq Y$ is a closed, convex cone,   $K^*$ is its dual cone, and 
\mbox{$s \leq_K t \,$}  means $t - s \in K. \,$ Note that the subscript 
$c$ refers to the objective  of the primal problem.

Problems $(\pc)$ and $(\dc)$ generalize linear programs and  share  some of the duality theory of linear
programming. For instance, a pair of feasible solutions  always satisfies the weak duality inequality 
$\la  c, x \ra \leq \la b, y \ra.$ 
However, in conic linear programming pathological phenomena occur: 
the optimal values of $(\pc)$ and of $(\dc)$ may differ, and they 
may not be attained. 

In particular, semidefinite programs (SDPs) and second order conic programs (SOCPs) 
--- probably the most useful and pervasive 
conic linear programs --- often behave pathologically:  for a variety of examples we refer to the textbooks \cite{BentalNem:01, Ren:01, BonnShap:00, Barvinok:2002, Tuncel:11}
surveys  \cite{VanBo:96, Todd:00, LuoSturmZhang:97} and research papers 
\cite{Ramana:97, AndRoosTerlaky:02, TunWolk:12}. 
Pathological conic LPs are both  theoretically interesting and often 
difficult, or impossible to solve numerically. 

These pathologies arise, since the  linear image of a closed 
convex cone is not always closed. 
For recent studies about when such sets are closed (or not), see e.g., 
\cite{BausBor:99, Aus:96, Pataki:07}. Three approaches (which we review in detail below) can help to avoid or  remedy the pathologies:  one can impose a constraint qualification (CQ), such as Slater's condition;
one can regularize $(\pc)\mhyphen(\dc)$ using a facial reduction algorithm \cite{BorWolk:81, WakiMura:12, Pataki:13}; or write an {\em extended dual}  \cite{Ramana:97, KlepSchw:13},  which uses 
extra variables and constraints. However, such CQs often do not hold, and neither 
facial reduction algorithms nor extended duals can help solve all pathological instances.

We started this research observing that pathological SDPs in the literature look curiously similar and 
one of our main goals is to find the root cause of the similarity.
We focus on the system underlying 
$(\pc)$ and call 
\beq \tag{$P$} \label{p}
\ba{rcl} 
\A x  \leq_K  b 
\ena
\eeq
{\em badly behaved} if there exists $c$ such that 
$(\dc)$ either does not attain its value or its value differs from the value 
of $(\pc).$ 
 We call (\ref{p}) {\em well behaved} if it is not badly behaved.

\nin{\bf Main contributions of the paper:} 
\benum

\item In Theorem \ref{unif-d-thm} of Section \ref{section-p} we characterize 
 when the system \eref{p} is
badly or well behaved. 
At the heart of Theorem \ref{unif-d-thm} is a simple geometric condition that involves the set of
feasible directions at $z \in K, \,$ i.e., 
$$  
\{ \, y \, | \, z + \eps y \in K \, \text{for some} \, \eps > 0 \, \}, 
$$ 
and $z$ is chosen as  a certain {\em slack}  in (\ref{p}). 

\co{Theorem \ref{unif-d-thm} implies in a surprisingly simple, and straightforward manner that two well-known conditions -- Slater's condition, and $K$ being polyhedral -- are sufficient for \eref{p} to be well behaved. }

In Theorem \ref{unif-d-thm} we unify two well-known (and seemingly unrelated) conditions for 
(\ref{p}) to be well behaved: the first is Slater's condition, and the second requires  $K$ to be polyhedral.

Theorem \ref{unif-d-thm} relies on a result on the closedness of the linear image of a closed convex cone from \cite{Pataki:07} (which we recap in Lemma \ref{pataki-cl}). 

		\item  In Section \ref{section-proof-sdp} we 
		 characterize when a semidefinite system 
	\beq \label{p-sd} \tag{\mbox{$P_{SD}$}}
	\sum_{i=1}^m x_i A_i \preceq B 
	\eeq
	is badly behaved via certain {\em excluded matrices}. 
	We assume (with no  loss of generality)
	that  
	a maximum rank positive semidefinite matrix of the form  $B - \sum_i x_i A_i$
	is 
	\beq \label{excluded}
	Z \, = \, \bpx I_r & 0 \\ 0 & 0 \epx  \, {\rm for \, some} \, 0 \leq r \leq n.
	\eeq
	We prove (in Theorem \ref{badsdp}) that \eref{p-sd} 
	is badly behaved iff there is a matrix $V$ which is a linear combination of 
	the $A_i$ and $B$ of the form 
	\beq \label{excluded2}
	V \, = \, \begin{pmatrix} V_{11} & V_{12} \\
		V_{12}^T & V_{22} 
	\end{pmatrix}, \, 
	\eeq
	where $V_{11}$ is $r \! \times \! r, \,$ $V_{22}$ is positive semidefinite, and $\R(V_{12}^T) \not \subseteq \R(V_{22}).$ Here 
	$\R()$ stands for rangespace. 
	
	The  excluded matrices $Z$ and $V$ 
	are easy to spot in all published badly behaved semidefinite systems (we counted about $20$  in the above references). 
	The simplest such system is 
	\beq \label{ex1-alpha}
	\ba{rl}
	x_1 \bpx \alpha & 1 \\ 1 & 0 \epx \preceq \bpx 1 & 0 \\ 0 & 0 \epx,
	\ena
	\eeq
	where $\alpha$ is any 
	real number: in (\ref{ex1-alpha}) 
	 the right hand side serves as $Z$ and the matrix on the left hand side serves as $V.$ 
	
	Theorem \ref{goodsdp} similarly characterizes well behaved semidefinite systems.
	
	Theorems \ref{badsdp} and \ref{goodsdp} follow from 
	Theorem \ref{unif-d-thm},  and Lemma \ref{psd-dir}, which characterizes the set of feasible directions and related sets in the semidefinite cone.

\item How do we verify that (\ref{p-sd}) is badly or well behaved? 
In other words, how do we convince a nonexpert 
reader that an instance of 
(\ref{p-sd}) is badly or well behaved? 
Theorems \ref{badsdp-re} and \ref{goodsdp-re} in Section \ref{npcapconp} show  
how to transform \eref{p-sd} into an equivalent standard system, 
whose bad or good  behavior is self-evident.
The transformation is surprisingly simple, as it 
 relies mostly on  elementary row operations --- the same operations that are used  in Gaussian elimination.
A natural analogy (and  our inspiration)  is how 
one transforms an infeasible 
linear system of equations $Ax = b$ to derive the obviously infeasible 
equation $\la 0, x \ra = 1.$

Here we also prove that i) 
badly/well behaved semidefinite systems are in 
$\np \cap \conp$ in the real number model of computing 
ii) for a well behaved semidefinite system we can restrict  optimal dual matrices
to be block-diagonal, and  iii) roughly speaking, we can partition a well behaved 
system into a strictly feasible part, and a linear part.

As a byproduct, we prove that all linear maps that act on symmetric matrices 
can be brought into a canonical form; this canonical form allows us
to easily check whether the image of the semidefinite cone under the given linear map is closed.  

\item In Section \ref{section-conclude} we  sketch analogous results for 
conic linear programs and SDPs in the dual form and prove  that all badly behaved semidefinite systems can be reduced,
by a sequence of natural operations, to the system \eref{ex1-alpha}. 

\item Since most examples in the main body of the paper have at most three variables and 
$3 \! \ti \! 3$ matrices, in Appendix A we give  a larger illustrative example with four
variables and $4 \! \ti \! 4$ matrices. We prove other technical results  in Appendix B. 

\eenum

 We illustrate our results by many examples.
 The only technical proofs in the main body of the paper are those of 
Theorem \ref{unif-d-thm} and of Lemma \ref{max-sl-lemma}, 
and these  can be safely skipped at first reading.

\nin{\bf Related work}  A fundamental question in convex analysis is
whether the linear image of a closed convex cone is closed. In this paper we rely on 
Theorem 1.1 from \cite{Pataki:07}, which we summmarize in Lemma \ref{pataki-cl}. This result 
gives several necessary conditions,  and 
exact characterizations for the class of {\em nice} cones. 
We refer to Bauschke and Borwein \cite{BausBor:99} for the closedness of the continuous 
image of a closed convex cone; to 
Auslender \cite{Aus:96} for the closedness of the linear image of an arbitrary closed convex set;
and  to Waksman and Epelman \cite{WaksEpel:76} for another related result.
For perturbation results we refer to 
Borwein and Moors \cite{BorweinMoors:09, BorweinMoors:10}; the latter paper shows that the 
set of linear maps under which the image of a closed convex cone is {\em not} closed  is 
small both in terms of measure and category. For a more general problem, 
whether the intersection of an infinite sequence of nested sets is nonempty, 
Bertsekas and Tseng \cite{BertTseng:07} gave a sufficient condition.
Their characterization is in terms of a certain {\em retractiveness} property of the set sequence.  

We say that $(\dc)$ is a {\em strong dual of} $(\pc)$ if they have the same value, and $(\dc)$ attains this value, when it is finite. Thus in general $(\dc)$ is not a strong dual of $(\pc).$ Using this terminology, 
(\ref{p}) is well behaved exactly if $(\dc)$ is a strong dual of $(\pc)$ for {\em all} $c.$  
We say that (\ref{p}) satisfies {\em Slater's condition}, if there is $x$ such that $b - \A x $ is in the relative interior of $K;$ if this condition holds, then (\ref{p}) is well behaved. 

Ramana in \cite{Ramana:97} proposed a strong dual for SDPs, 
which uses polynomially many extra variables and constraints.
His result implies that semidefinite feasibility 
is in $\np \cap \conp$ in the real number model of computing.
Klep and Schweighofer in 
\cite{KlepSchw:13} constructed 
a Ramana-type strong dual for SDPs, 
which, interestingly, is 
based on ideas from algebraic geometry, rather than from  convex analysis. 

The facial reduction algorithm of Borwein and Wolkowicz in \cite{BorWolk:81,BorWolk:81B} 
converts \eref{p} into a system that satisfies Slater's condition, and is hence well behaved.
The algorithm relies on a sequence of reduction steps.
For more recent, simplified facial reduction algorithms, see
Waki and Muramatsu \cite{WakiMura:12} and Pataki \cite{Pataki:13}. 
Ramana, Tun\c{c}el, and Wolkowicz in \cite{RaTuWo:97} 
proved the correctness of Ramana's dual from 
the facial reduction algorithm of 
\cite{BorWolk:81, BorWolk:81B}, showing the connection of these two seemingly unrelated 
concepts. We refer to Ramana and Freund \cite{RaFreund:96} 
for a proof that the  Lagrange dual of Ramana's dual has the same value as the original problem. Generalizations 
of Ramana's dual are known for conic LPs over {\em nice} cones \cite{Pataki:13};  
and for conic LPs over homogeneous cones (P\'olik and Terlaky \cite{PolikTerlaky:09}).

For a generalization of the concept of strict 
complementarity  (a concept that plays an important role in our work),
we refer to Pena and Roshchina \cite{PenaRosh:13}. 
Schurr et al in \cite{Schurretal:07} characterize  
{\em universal duality} ---  when strong duality holds
for all right hand sides and objective functions. 
Tun\c{c}el and  Wolkowicz in \cite{TunWolk:12} related the lack of strict complementarity 
in a homogeneous conic linear system to the existence of an objective 
function with a positive gap. 

We finally remark that the technique of {\em reformulating} equality constrained
SDPs (relying mostly on elementary row operations),  to easily verify their infeasibility 
was used recently by Liu and Pataki \cite{LiuPataki:15}.

\subsection{Preliminaries. When is the linear image of a closed convex cone closed?} 
We now review some basics in convex analysis, relying mainly on 
references \cite{Rockafellar:70, HirLemar:93, BorLewis:00, BausComb:11}. 
In Lemma \ref{pataki-cl} we also give a short and transparent 
summary of a result on the closedness 
of the linear image of a closed convex cone from \cite{Pataki:07}. 

If $x$ and $y$ are elements of the same Euclidean space, we sometimes 
write $x^*y$ for $\la x, y \ra.$ 
For a set $C$ we denote its linear span, the orthogonal complement of 
its linear span, its closure, and interior by 
$\lin C, \, C^\perp, \, \cl C, \,$ and $\myint C, \,$ respectively.
For a convex set $C$ we denote its relative interior by 
$\ri C.$ 
For a convex set $C \,$ and $x \in C \,$ we define 
\beqa \label{dirldir}
\dir(x, C) & = & \{ \, y \, | \, x + \eps y \in C \, \text{for some} \, \eps > 0 \, \}, \\ \label{defldir}
\ldir(x, C) & = & \dir(x, C) \cap - \dir(x, C), \\ \label{deftan}
\tan(x, C)  & = & \cl \dir(x, C) \cap - \cl \dir(x, C).
\eeqa
Here $\dir(x, C)$ is the set of feasible directions at $x$ in $C, \,$  
and $\tan(x, C)$ is the tangent space at $x$ in $C.$

A set $C$ is a {\em cone } if $\lambda x \in C$ holds for all $x \in C, \,$ and $\lambda \geq 0.$ 
Let $C$ be a  closed convex cone. Its dual cone is 
$$
C^* \, = \, \{ \, y \, | \, \la y, x \ra \geq 0 \,\, \forall x \in C \}.
$$
For $E, \,$  a convex subset of $C, \,$ we say that 
$E$ is a {\em face} of $C, \,$ 
if $x_1, x_2 \in C, \,$ and $1/2(x_1 + x_2) \in E \,$ implies that 
$x_1$ and $ x_2$ are in $E.$ 
\co{An equivalent definition is that 
$x_1, \, x_2 \in C, \, x_1 + x_2 \in E$ implies that 
$x_1$ and $ x_2$ are in $E.$ } 
\co{For an $E$ face of $C$ we define 
its {\em conjugate face} 
as 
\beqa \nonumber
E^\triangle & = & C^* \cap E^\perp,
\eeqa
and the conjugate face of a $G$ face of $C^*$ as $G^\triangle \, = \, C \cap G^\perp.$
}

For $\, x \in C \,$ and $u \in C^*, \, $ we say that $u$ is 
{\em strictly complementary to} $\, x \,$ if 
$u \in \ri (C^* \cap x^\perp).$ 
If $C$ is the semidefinite cone, or the second order cone, then  
 $u$ is strictly complementary to $x \,$ iff $x$ is 
  strictly complementary to $u;$ in other cones, however, this may not be the case
  (see a discussion  in \cite{Pataki:00A}).

We say that a closed convex cone $C$ is {\em nice}, if 
\beq
\nonumber 
C^* + E^\perp   \; {\rm is \; closed } \; {\rm for \; all \; } E \; {\rm faces \; of \;} C.
\eeq

We know that polyhedral, semidefinite, and $p$-order cones are nice \cite{BorWolk:81, BorWolk:81B, Pataki:07}; 
the intersection of a nice cone with a linear subspace 
and the linear preimage of a nice cone are nice \cite{ChuaTuncel:08}; hence 
homogeneous cones are nice, as they are the 
intersection of a semidefinite cone with a linear subspace 
(see \cite{Chua:03, Faybu:02}). 
In \cite{Pataki:12} we characterized 
nice cones, proved that they must be facially exposed and conjectured that all facially exposed cones 
are nice. However, Roshchina \cite{Vera:13} disproved this conjecture. 

We denote the rangespace, nullspace, and adjoint operator of a 
linear operator $\M$ by $\R(\M), \, \N(\M)$ and $\M^*, \,$ respectively. 
We denote by $\sym{n}$ the set of $n$ by $n$ symmetric matrices, and 
by $\psd{n}$ the set of $n \times n$ symmetric positive semidefinite (psd) matrices.
For symmetric matrices $A$ and $B$ we write $A \preceq B$ [$A \prec B$] 
to denote that 
$B - A$ is positive semidefinite [positive definite], 
and we write $A \bullet B$ to denote the trace of 
$AB.$ We have $(\psd{n})^* = \psd{n}$ with respect to the $\bullet$ inner product.

We will use the  fact that for an $H \subseteq \sym{n}$ affine subspace  
$$
\ri (H \cap \psd{n}) \, = \, \{ \, X  \in \psd{n} \, | \, X \, {\rm is \, a \, maximum \, rank \, psd \, matrix \, in} \, H \}.
$$

For $A, B \in \sym{n}$ and an  invertible matrix $T$  we will use the identity 
\beq
\label{txt}
T^T A T \bullet T^{-1} B T^{-T} \, = \, A \bullet B.
\eeq
For matrices $A_1$ and $A_2, \,$ we let 
$$
A_1 \oplus A_2 \, = \, \bpx A_1 & 0 \\
0  & A_2 \epx,
$$
and for sets of matrices $X_1$ and $X_2$ we define 
$$
X_1 \oplus X_2 \, = \, \{ \, A_1 \oplus A_2 \, | \, A_1 \in X_1, \, A_2 \in X_2 \, \}.
$$
For instance, $\psd{r} \oplus \{ 0 \}$ (where the order of the $0$ matrix will be clear from 
context) is the set of matrices with the upper left 
$r \times r$ block positive semidefinite and the rest of the components zero. 

We write $I_r$ for the identity matrix of order $r.$

The following question is fundamental in convex analysis: when is 
the linear image of a closed convex cone closed? 
We state and illustrate a short version of Theorem 1.1 from \cite{Pataki:07}, which gives easily checkable
conditions which are ``almost''  necessary and sufficient. We will use Lemma \ref{pataki-cl} later on to prove 
Theorem \ref{unif-d-thm}. 
\ble \label{pataki-cl}
Let $\M$ be a linear map, $C$ a closed convex cone, and $w \in \ri ( C \cap \R(\M) ).$
Conditions \eref{cl2} and \eref{cl3} below 
are equivalent to each other, and {\em necessary} 
for $\M^*C^*$ to be closed. If $C$ is nice, then they are {\em necessary and sufficient.} 
\benum
\item \label{cl2} $\R(\M) \cap \bigl(\cl \dir(w, C) \setminus \dir(w,C)\bigr) = \emptyset.$ 
\item \label{cl3} There is $w' \in \N(\M^*) \cap C^* \,$ strictly complementary to $w, \,$ and 
$$\R(\M) \cap \bigl(\tan(w,C) \setminus \ldir(w,C)\bigr) = \emp.$$
\eenum
\ele
\qed

Our first example which illustrates Lemma \ref{pataki-cl} is very simple: 
\bex \label{1stex} {\rm 
	 Let $C = C^* = \psd{2}$ and 
	define the map $\M: \rad{2} \rightarrow \sym{2}$ as
	$$
	\M(x_1, x_2) \, = \, \bpx  x_1 & x_2 \\ x_2 & 0 \epx.
	$$
	Then 
	$\M^*Y = (y_{11}, 2 y_{12})^T \, {\rm where} \, Y \in \sym{2}, \,$ and $\M^* C^*$ is not closed:
	a direct computation shows 
	 $\M^*C^* = (\rad{}_{++} \ti \rad{}) \cup (0,0),$ where 
	$\rad{}_{++}$ stands for the set of strictly positive reals. 
	
	Lemma \ref{pataki-cl} also proves that  $\M^*C^*$ is not closed: to see how, let 
	$$
	w = \bpx 1 & 0 \\ 0 & 0 \epx, \, v = \bpx 0 & 1 \\ 1 & 0 \epx.
	$$
	Then $w \in \ri(\R(\M) \cap C), \,$ since it is a maximum rank psd matrix in $\R(M). \,$	Also, 
	$v \in \R(\M) \cap (\cl \dir(w,C) \setminus \dir(w,C)), \,$ since 
	$v \not \in \dir(w,C)$ follows from the definition, and 
	$v \in \cl \dir(w,C)$ follows, since putting any $\eps >0$ into the 
	$(2,2)$ position of $v$ makes it a feasible direction. So condition  (\ref{cl2}) of Lemma 
	\ref{pataki-cl} is violated, hence $\M^*C^*$ is not closed. } 
	\eex
		The next, more involved example illustrates the key point of Lemma 
	\ref{pataki-cl}: the image set $\M^*C^*$ usually has much 
	more complicated geometry than   $C$ and  $C^*.$ 
Lemma \ref{pataki-cl} sheds light on the geometry of $\M^*C^*$ via the geometry of the simpler set  $C.$ 
\bex{\rm \label{MCex} 
	Let $C = C^* = \psd{3}$ and 
	define the map $\M: \rad{3} \rightarrow \sym{3}$ as 
	$$
	\M(x_1, x_2,x_3) = \bpx x_1 & 2 x_2 & x_3 \\ 2 x_2 & x_2 + x_3 & 0 \\ x_3 & 0  & 0 \epx. 
	$$
Thus $\M^*Y = (y_{11}, y_{22}+ 4 y_{12}, y_{22}+2y_{13}), \, $ where $Y \in \sym{3}.$

	It is a straightforward computation (which we omit) to show 
	\beq \label{MsCs} 
	\ba{rcl} 
	\cl (\M^* C^*) & = & \{ (\alpha, \beta, \gamma): \, \alpha \geq 0, 4 \alpha + \beta \geq 0 \} ,\\
		\cl (\M^* C^*) \setminus \M^* C^* & = & \{ (0, \beta, \gamma): \, \gamma \neq \beta \geq 0 \}. 
		\ena
		\eeq
		The set $\M^* C^*$ is shown on Figure \ref{frontier} in blue, 
		and ${\rm cl \,} (\M^*C^*) \setminus \M^*C^*$ 
	 in green. (Note that the blue diagonal segment on the green facet actually belongs to $\M^*C^*.$)
		
	Lemma \ref{pataki-cl} easily proves 
		that  $\M^*C^*$ is not closed, even without computing the sets in (\ref{MsCs}); indeed, let 
$$
	w \, := \, \M(6,1,0) \, = \, \bpx 6 & 2 & 0 \\ 2 & 1 & 0 \\ 0 & 0 & 0 \epx,  \, \, v := \M(0,0,1) \, = \, \bpx 0 & 0 & 1 \\ 0 & 1 & 0 \\ 1 & 0 & 0 \epx,
	$$
	and observe i) 	$w \in \ri (\R(\M) \cap C),$ since it is a maximum rank psd matrix in $\R(\M);$   ii) $v \not \in \dir(w,C)$ follows from the definition; 
	and  iii) $v \in \cl \dir(w,C),$ since putting any $\epsilon > 0$ into the $(3,3)$ position of $v$ makes it a feasible direction. 
	
	Thus condition \eref{cl2} in Lemma \ref{pataki-cl} is violated, so 
	$\M^*C^*$ is not closed. 
}
	\eex

\begin{figure}[H]
	\centering
	\includegraphics[scale=0.4]{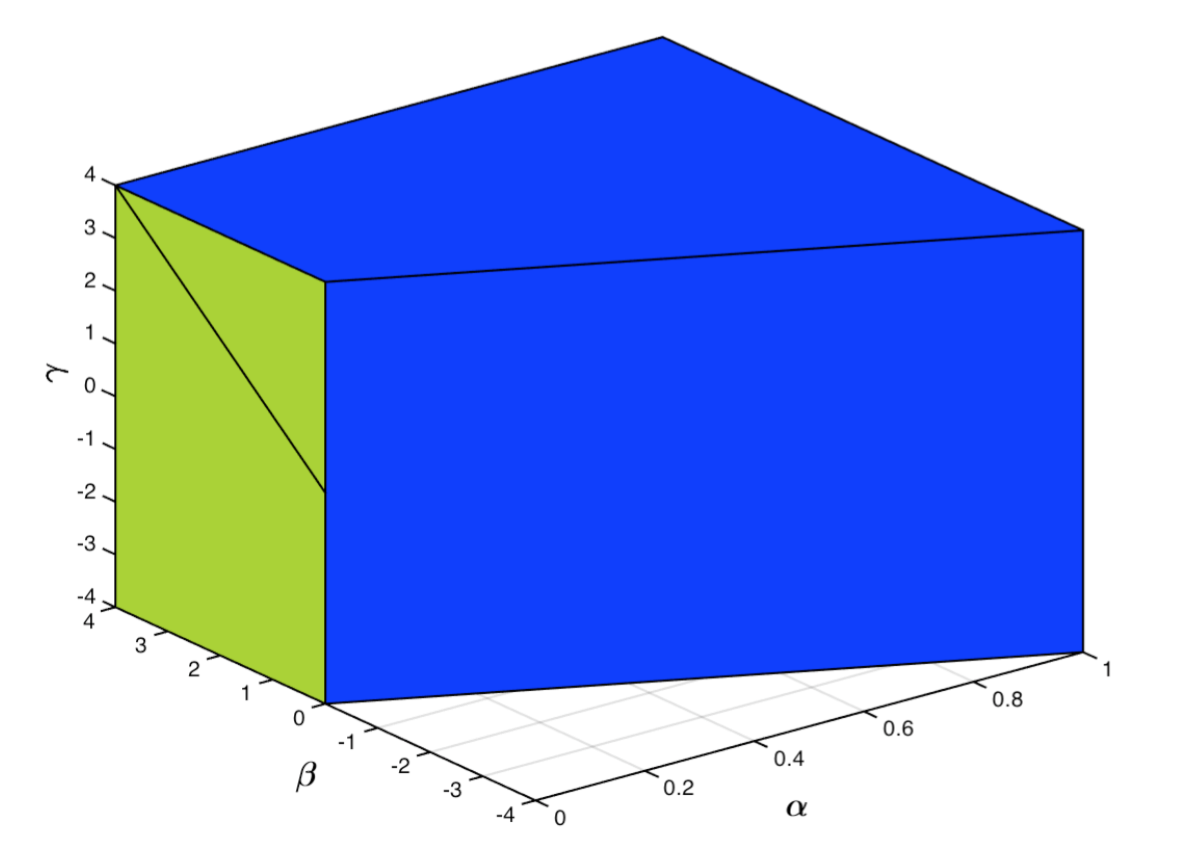}
	\caption{The set $\M^*C^*$ is in blue, and 
		${\rm cl} \, (\M^* C^*) \setminus \M^*C^*$   is in green} 
	\label{frontier}
\end{figure}

\co{
\begin{figure}[H]
	\centering
	\includegraphics[scale=0.4]{frontier9}
	\caption{The set $\M^*C^*$ is in blue, and 
		${\rm cl} \, (\M^* C^*) \setminus \M^*C^*$   is in green} 
	\label{frontier}
\end{figure}

\begin{figure}[H]
	\centering
	\includegraphics[scale=0.4]{frontier16}
	\caption{The set $\M^*C^*$ is in blue, and 
		${\rm cl} \, (\M^* C^*) \setminus \M^*C^*$   is in green} 
	\label{frontier}
\end{figure}
}

We mention in passing that the second part of condition \eref{cl3} in Lemma \ref{pataki-cl} 
is stated in Theorem 1.1 in \cite{Pataki:07} as
$\R(\M) \cap ((E^{\triangle})^\perp \setminus \lin E) = \emp,$ where 
$E$ is the smallest face of $C$ that contains $w$ and $E^\triangle = C^* \cap w^\perp.$
However, this is an equivalent formulation, as implied by the 
characterization of $\lin E$ and $(E^{\triangle})^\perp$,  see e.g., Lemma \ref{ldir} in Appendix B. 
\co{Second, Theorem 1.1 in \cite{Pataki:07} actually shows that 
 conditions \eref{cl2} and \eref{cl3} suffice  for 
$\M^*C^*$ to be closed under the weaker assumption that $C^* + E^\perp$ is closed, where 
$E$ is the smallest face of $C$ that contains $w.$ Hence Corollary \ref{cor-\MC} follows 
from Theorem 1.1 in \cite{Pataki:07} even if we do not assume that $C$ is nice.
}


\co{For optimization problems we use the symbol $\val()$ to denote their optimal value. 
For program $(\dc)$ we define 
$\{ y_i \} \subseteq K^*$ to be an {\em asymptotically feasible (AF)} solution,
if $ A^* y_i \rightarrow c, $ and the {\em asymptotic value of $(\dc)$} is 
$$
\ba{rcl}
\aval(\dc) & = & \inf \{ \, \lim \, b^* y_i \, | \, \{ y_i \} \, {\rm is \; asymptotically \; feasible \; to \;\;} (\dc) \, \},
\ena
$$ 
where the infimum is taken over those AF solutions for which 
$\lim \, b^* y_i$ exists. 
\ble \label{asd} (Duffin \cite{Duff:56}) 
Problem $(\pc)$ is feasible with $\val(\pc) < + \infty, \,$ iff 
$(\dc)$ is asymptotically feasible with $\aval(\dc) > - \infty, \,$ and if these equivalent statements hold,  
then 
\beq \nonumber 
\val(\pc) = \aval(\dc).
\eeq
\qed
\ele
}

Throughout the paper we assume that \eref{p} is feasible.
Recall that we  say that \eref{p} satisfies {\em Slater's condition} if there exists 
$x$ such that $b - \A x \in \ri K.$

\section{When is a conic linear system  badly or well behaved?} %

\label{section-p} 

In this section we present our main characterization of when \eref{p} is badly or well behaved (these concepts are defined in the Introduction).
We first need a definition.
\begin{Definition} \label{maxslack}
A {\em  slack} in \eref{p} is a vector in 
$$
(\R(\A) + b ) \cap K, 
$$
and a {\em maximum slack} is a vector in the relative interior of all slacks.
\end{Definition}

We start with a basic lemma:
\co{ that connects the good behavior of \eref{p} with the closedness 
of $K^* \ti \rad{}_+$ under a certain linear map. }
\ble \label{duffin}
The system \eref{p} is well behaved, if and only if 
the set
$$
\left(\!\!\! \begin{array}{cc}
   \A^* & 0 \\
   b^* & 1
\end{array} \!\!\! \right)  \left( \!\!\! \begin{array}{cc}
   K^* \\
    \mathbb{R}_+
\end{array}  \!\!\! \right)
$$
is closed.
\ele
\qed

To put Lemma \ref{duffin} into perspective, note that  
the image set in Lemma \ref{duffin} is closed  if $(\A,b)^*K^*$ is closed
(one can argue  this directly  or by modifying the proof of Lemma \ref{duffin}).
In turn, if  $(\A,b)^*K^*$ is closed, then the duality gap between $(\pc)$ and $(\dc)$ is zero, even if $K$ lives in an infinite  dimensional  space ---
see, e.g., Theorem 7.2 in \cite{Barvinok:2002} (where our primal is called the dual).
The proof of Lemma \ref{duffin} is  standard, and we give it in Appendix B.

The main result of this section follows (recall the definition of $\dir(z,K)$ and related sets from 
\eref{dirldir}--\eref{deftan}). We  write $\R({\cal A}, b)$ for the 
rangespace of the operator
$(x,t) \rightarrow {\cal A}x + bt.$

\bth \label{unif-d-thm} Let $z$ be a maximum slack in \eref{p}. 
Conditions \eref{frdir} and \eref{tang} below 
are equivalent to each other, and {\em necessary} 
for \eref{p} to be well behaved.  If $K$ is nice, then they are {\em necessary and sufficient.} 
\benum 
\item \label{frdir} $\R({\cal A}, b) \cap (\cl \dir(z, K) \setminus \dir(z,K))  = \emp. \,$ 
\item \label{tang}  
There is $u \in \N\bigl((\A,b)^*\bigr) \cap K^* \,$ strictly complementary to $z, \,$ and 
$$\R(\A,b) \cap \bigl(\tan(z,K) \setminus \ldir(z,K)\bigr) = \emptyset.$$
\eenum
\enth
\qed

To build intuition we  show 
how Theorem \ref{unif-d-thm} unifies two classical, seemingly unrelated, 
{\em sufficient} conditions for \eref{p} to be well-behaved.
\bcor \label{well-cor} 
Suppose that $K$ is a nice cone. If $K$ is polyhedral or (\ref{p}) satisfies Slater's condition, 
then \eref{p} is well-behaved.
\ecor
\pf{} Let $z$ be a maximum slack in (\ref{p}). 
If $K$ is polyhedral, then so is  $\dir(z,K)$. If (\ref{p}) satisfies Slater's condition, then clearly 
$z \in \ri K,$ so $\dir(z,K) = \lin K. \,$ In both  cases $\dir(z,K)$ is closed, hence 
Condition  \eref{frdir} holds, so  \eref{p} is well behaved. \qed

\co{
\pf{} First suppose that $K$ is polyhedral. Then $\dir(z,K)$ is closed for any $z \in K, $  hence Condition \eref{frdir} in Theorem \ref{unif-d-thm} holds, so \eref{p} is well behaved.  
Next assume $z \in \ri K.$ Then 
 $\dir(z,K) = \lin K, \,$ which is 
a closed set, so Condition  \eref{frdir} again holds, therefore \eref{p} is well behaved. \qed
}

Though Lemma \ref{duffin} is a bit simpler to state than Theorem \ref{unif-d-thm},
the latter  will be more useful.
On the one hand,  Lemma \ref{duffin} relies on the closedness of the linear image of $K^* \times \rad{}_+, \,$ 
which may not be easy to check.  
On the other hand, Theorem \ref{unif-d-thm} relies on the geometry of the cone $K \,$ itself, and not on the geometry of its linear image. The geometry of typical cones that occur in optimization
--- e.g. the geometry of  the semidefinite cone --- is well understood.
Thus Theorem \ref{unif-d-thm}, among other things, will lead to a proof 
that badly behaved semidefinite systems  are 
in $\np \cap \conp$ in the real number model of computing. 
Lemma \ref{duffin}, by itself, affords no such  corollary. 
%

Note that if $b=0$, then by Lemma \ref{duffin} the system \eref{p} is well behaved  iff 
$\A^*K^*$ is closed. Thus in this case Lemma \ref{pataki-cl} and Theorem \ref{unif-d-thm} 
are equivalent. To prove the general case of Theorem \ref{unif-d-thm} 
we  use a homogenization argument.

\co{
To prove Theorem \ref{unif-d-thm} we will combine Lemmas \ref{pataki-cl} and \ref{duffin}.
Note that if $b=0$, then by Lemma \ref{duffin}, the system \eref{p} is well behaved, iff 
$\A^*K^*$ is closed. Thus in this case Lemma \ref{pataki-cl} 
}

\co{
the  then Lemma \ref{pataki-cl} and Theorem \ref{unif-d-thm} are equivalent: in this case an elementary argument shows that 
$(\pc)$ is bounded, iff $c \in \cl \A^*K^*.$ Hence in the  $b=0$ case 
\eref{p} is well behaved iff $\A^*K^*$ is closed. To prove the general case of Theorem \ref{unif-d-thm} 
we need a homogenization 
argument.
}

\pf{of Theorem \ref{unif-d-thm}:} 
We consider the homogenized system 
\beq \tag{$P_{h}$} \label{phomplus}
\ba{rcl} 
\A x - b x_0 & \leq_K & 0 \\
- x_0        & \leq & 0, 
\ena
\eeq
and first prove the following claim: 

{\bf Claim} There is a $z$  maximum slack in \eref{p} such that $(z,1)$ is also a maximum slack in \eref{phomplus}.

To prove the Claim we first note that if $z$ is a maximum slack in \eref{p} and 
$z'$ is some other slack, then $\lambda z  + (1-\lambda)z'$ is also a maximum slack for all 
$0 < \lambda \leq 1 \,$ (by Theorem 6.1 in \cite{Rockafellar:70}).
A similar result holds for \eref{phomplus}.

Now let $z_1$ be a maximum slack in \eref{p}, then $(z_1,1)$ is a slack in 
\eref{phomplus}. Next, let $(z_2, x_0)$ be a maximum slack in \eref{phomplus}.
By the properties of the relative interior, and since $(z_1, 1)$ is a slack in 
\eref{phomplus}, we have that 
$(z_2, x_0) - \eps (z_1,1)$ is a slack in (\ref{phomplus}) for some $\eps >0. \,$ 
So $x_0 > 0$ must hold, and (after normalizing) we can 
assume $x_0=1.$ Hence $z_2$ is a slack in \eref{p} and 
$$
z := \dfrac{1}{2} (z_1 + z_2)
$$
will do. This completes the proof of the claim.

To proceed with 
the proof of the theorem, we note that the set of maximum slacks in \eref{p} is a relatively open set, so by Theorem 
18.2 in \cite{Rockafellar:70} it is contained in $\ri F, \,$ where $F$ is some face of $K.$
Therefore   $\dir(z,K) = K + \lin F$ for any maximum slack $z$ (see e.g. Lemma 2.7 in \cite{Pataki:00A}) so the sets $\dir(z,K)$ and $\tan(z,K)$ depend only on $F.$ Hence 
we are free to use any maximum slack of \eref{p} in our proof, and we will use the particular 
maximum slack provided in the preceding Claim.

For convenience we define the linear map
$$
\A_{h} = \left(\!\!\! 
\begin{array}{cc}
   \A & b \\
   0 & 1
\end{array} \!\!\! \right) \, 
$$
which corresponds to the homogenized conic linear system \eref{phomplus}.

We first note that (trivially) 
\beqa \nonumber 
\dir \bigl( (z,1), K \times \rad{}_+ \bigr) & = & \dir(z, K) \times \rad{} \, {\rm holds}. 
\eeqa
Equations \eref{dirldir}--\eref{deftan} imply that the same 
statement holds, if we replace the operator $''\dir''$ by 
$''\cl \dir'', \,$ $ ``\tan'', \,$ or $''\ldir''. $ 

Hence the following equations hold: 
\beqa \label{h2.5} 
\cl \dir \bigl( (z,1), K \times \rad{}_+ \bigr)  \setminus \dir \bigl( (z,1), K \times \rad{}_+ \bigr) & = & \bigl( \cl \dir(z, K) \setminus \dir(z, K)\bigr) \times \rad{}, \\ 
\label{h4.5} 
\tan \bigl( (z,1), K \times \rad{}_+ \bigr) \setminus \ldir \bigl( (z,1), K \times \rad{}_+ \bigr) & = & \bigl( \tan(z, K) \setminus \ldir(z, K) \bigr) \times \rad{}. 
\eeqa
Consider now the following variants of conditions \eref{frdir} and \eref{tang}:
\benum
\item[$(1')$] \label{hom2} $\R(\A_{h}) \cap \bigl[ \cl \dir \bigl((z,1), K \times \rad{}_+ \bigr) \setminus \dir \bigl((z,1), K \times \rad{}_+\bigr)  \bigr] = \emptyset.$
\item[$(2')$] \label{hom3} There is $(u, u_0) \in \N(A_{h}^*) \cap (K \times \rad{}_+)^* \,$ strictly complementary to 
$(z,1)$ and 
$$
\R(\A_{h}) \cap \bigl[ \tan \bigl((z,1), K \times \rad{}_+ \bigr) \setminus \ldir \bigl((z,1), K \times \rad{}_+ \bigr)\bigr]  = \emptyset.
$$
\eenum
Since $(z,1)$ is a maximum slack in (\ref{phomplus}), we have  
$(z,1) \in \ri \bigl(\R(\A_h) \cap (K \times \rad{}_+) \bigr).$ 
Hence by Lemma  \ref{pataki-cl} with 
$
C = K \ti \rad{}_+, \, \M = \A_{h}, \, w = (z,1)
$ 
we find 
$$
\A_h^* (K \times \rad{}_+)^* \,\, {\rm is \,\, closed \,\,}  \RA (1') \LRA (2')
$$
and that equivalence holds when $K \ti \rad{}_+$ is nice. 

We next note that by \eref{h2.5}  condition  
$(1')$ 
is equivalent to \eref{frdir}. Also, if $(u, u_0)$ is as specified 
in $(2'),$ then 
$$
\la (u, u_0), (z,1) \ra \, = \, \la u,z \ra + u_0 \, = \, 0,
$$
and since both terms above are nonnegative, we must have $u_0 = 0. \,$ Thus using \eref{h4.5} we find that statement $(2')$ is equivalent to condition \eref{tang} in Theorem \ref{unif-d-thm}.  
Thus we have 
$$
\A_h^* (K \times \rad{}_+)^* \,\, {\rm is \,\, closed \,\,}  \RA (1)  \LRA (2) 
$$
with equivalence holding when $K \ti \rad{}_+$ is nice. 
Finally, $K$ is nice if and only if $K \times \rad{}_{+}$ is, thus invoking
Lemma \ref{duffin}  completes the proof. 
\qed

We  can easily modify the proof of  Theorem \ref{unif-d-thm} to show that 
conditions \ref{frdir} and \ref{tang} suffice for \eref{p} to be well behaved, even under a weaker condition than $K$ being nice: it
is enough for $K^* + F^\perp$ to be closed, where 
$F$ is the smallest face of $K$ that contains $z.$ This more general version of
Theorem \ref{unif-d-thm} implies that Corollary \ref{well-cor}
 holds even if we do not assume that  $K$ is nice --  we refer the interested reader to version 3 of the paper on arxiv.org.

\section{When is a semidefinite system  badly or well behaved?} 
\label{section-proof-sdp} 
We now specialize the results of Section \ref{section-p}, and characterize 
when the semidefinite system \eref{p-sd} is badly or well behaved. To this end, we consider the 
primal-dual pair of SDPs
\begin{center}                                                                                              
$$                                                                                                          
\begin{array}{lrlcrlr}                                                                                      
       &   \sup  & \sum_{i=1}^m c_i x_i            &    \hspace{2cm} & \inf   &   B \bullet Y  &   \\ 
(\sdpc)  &   s.t.  & \sum_{i=1}^m x_i A_i \preceq B  &    \hspace{2cm} & s.t.  &   Y \succeq 0   & (\sddc)     \\ 
       &         &                                 &    \hspace{2cm} &       &   A_i \bullet Y  = c_i \,\, (i=1, \dots, m), &    \end{array}                                                                                              
$$   
\end{center}                                                                                                      
where $A_1, \dots, A_m, B \in \sym{n},$ and $c_1, \dots, c_m$ are scalars.

Specializing Definition \ref{maxslack} to the semidefinite system 
(\ref{p-sd}), we find  that 
i) a  slack in (\ref{p-sd}) is a matrix of the form 
$ S = B - \sum_i x_i A_i \succeq 0, \,$  and ii) a maximum slack in (\ref{p-sd}) is a 
maximum {\em rank} slack. We also note that the cone of positive semidefinite matrices is nice \cite{BorWolk:81, BorWolk:81B, Pataki:07}.

We make the following 
\bass \label{slack-ass}
The maximum rank slack in \eref{p-sd} is 
\beq \label{Zslack}
Z \, = \, \bpx I_r & 0 \\ 0 & 0 \epx \, {\rm for \, some \,} 0 \leq r \leq n. 
\eeq
\eass

We can easily satisfy Assumption \ref{slack-ass}, at least from a theoretical point of view, 
as follows.  
If  $Z$ is any maximum rank slack in \eref{p-sd},  
$Q$ is a matrix of suitably scaled eigenvectors of $Z,$ and  we apply the rotation 
$Q^T()Q$ to all  $A_i$ and $B,$ then the  maximum rank slack in  
 the rotated system 
is in the required form. (We do not make a claim about actually computing $Z$ or $Q$; we discuss this point  more  at the end of Section \ref{npcapconp} ).

In the interest of the reader we first state and illustrate the main 
results, then prove them.
\bth \label{badsdp}
The system \eref{p-sd} 
is badly behaved 
if and only if there is a matrix $V$ which is a linear combination of the 
$A_i$ and $B$ of the form 
\beq \label{Vform}
V \, = \, 
\begin{pmatrix} V_{11} & V_{12} \\
                V_{12}^T & V_{22} 
\end{pmatrix}, \, 
\eeq
where $V_{11}$ is $r \times r, \,$ $V_{22} \succeq 0, \,$ and 
$\R(V_{12}^T) \not \subseteq \R(V_{22}).$ 
\qed \enth
The $Z$ and $V$ matrices provide a {\em certificate} of the bad behavior of 
\eref{p-sd}. 

\bex \label{ex1} {\rm In the problem 
\beq \label{ex1-problem}
\ba{rl}
\sup &  x_1 \\
s.t. & x_1 \bpx 0 & 1 \\ 1 & 0 \epx \preceq \bpx 1 & 0 \\ 0 & 0 \epx
\ena
\eeq
the only feasible solution is $x_1 = 0.$ 
The dual program, in which  we denote the components of $Y$ by $y_{ij}, \,$ 
is equivalent to 
\beq  \nonumber
\ba{rl}
\inf &  y_{11}  \\
s.t. & \bpx y_{11} & 1/2 \\ 1/2 & y_{22} \epx \succeq 0,
\ena
\eeq
which has a $0$ infimum but does not attain it. 

The certificates of the bad behavior of the system in \eref{ex1-problem} are
$$
Z \, = \, \bpx 1 & 0 \\ 0 & 0 \epx, \, V \, = \, \bpx 0 & 1 \\ 1 & 0 \epx.
$$
}
\eex

\bex \label{ex2}
{\rm The problem 
\beq \label{ex2-problem} 
\ba{rl}
\sup  &  x_2   \\
s.t. & x_1 \bpx 1 & 0 & 0 \\ 0 & 0 & 0 \\ 0 & 0 & 0 \epx + x_2 \bpx 0 & 0 & 1 \\ 0 & 1 & 0 \\ 1 & 0 & 0 \epx \preceq 
\bpx 1 & 0 & 0 \\ 0 & 1 & 0 \\ 0 & 0 & 0 \epx 
\ena
\eeq
again has an attained $0$ supremum. The reader can easily check that 
the value of the dual program is $1 \,$ (and it is attained), so 
there is a finite, positive duality gap. 

In \eref{ex2-problem} the right hand side is the maximum slack, 
and we can choose the coefficient matrix of $x_2$ 
as the $V$ matrix of Theorem \ref{badsdp}. 
}
\eex

We next characterize well behaved semidefinite systems: 
\bth \label{goodsdp}
The system \eref{p-sd} is well behaved 
if and only if conditions \eref{cond1} and \eref{cond2} below
hold: 
\benum
\item \label{cond1} There is a matrix $U$ of the form 
\beq \label{U} 
U \, = \, \bpx 0 & 0 \\ 0 & U_{22} \epx,
\eeq
with $U_{22} \in \sym{n-r}, \, U_{22} \succ 0$ and 
\beq \label{AiU} 
A_1 \bullet U \, = \, \dots \, = \, A_m \bullet U \, = \, B \bullet U \, = \, 0.
\eeq
\item \label{cond2} For all $V$ matrices, which are a linear combination of the $A_i$ and $B$ and are of the form 
\beq \nonumber
V \, = \, 
\bpx V_{11}   & V_{12} \\
     V_{12}^T & 0  \epx,
\eeq
with $V_{11} \in\sym{r},$ we must have  $V_{12} = 0.$ 
\eenum \qed
\enth

\bex \label{ex3} 
{\rm The system
\beq \label{seinfeld} 
x_1 \bpx 0 & 0 & 0 \\ 0 &  0 & 1 \\ 0 & 1 & 0 \epx \preceq \bpx 1 & 0 & 0 \\ 0 & 0 & 0 \\ 0 & 0 & 0 \epx
\eeq
is well behaved; we can easily prove this either directly or via Theorem \ref{goodsdp}.
To do the latter,  note that the right hand side of \eref{seinfeld} is the maximum rank slack,  condition \eref{cond1} of Theorem \ref{goodsdp} holds with 
$U = 0 \oplus I_2,$ 
and condition \eref{cond2} holds vacuously 
(the $(1,2), (1,3)$ block of both constraint matrices is zero). 
}
\eex

\bex
\label{alpha-ex}
{\rm This example
illustrates both badly and well  behaved semidefinite systems, depending on the value of the parameter
$\alpha:$ 
\begin{small}
\beq \label{lastsystem}
\ba{rcl}
x_1 \bpx 0 & 0 & 1   \\ 
              0 & 1 & -3 \\
              1 & -3 & 8 \epx 
+  x_2 \bpx 0 & 1 & -3  \\ 
             1 & 0 & 1 \\
             -3 & 1 & -6 \epx 
+ x_3 \bpx   1 & 1 & \alpha  - 3\\ 
             1 & 1 & -2 \\
   \alpha - 3 & -2 & 2 \epx 
& \preceq & \bpx 2 & 2 & \alpha - 5 \\ 2 & 2 & -4 \\  \alpha - 5 & -4 & 4 \epx
\ena
\eeq
\end{small}
Let us write $A_i$ for the constraint matrices on the left, and $B$ for the right hand 
side matrix in (\ref{lastsystem}). We first observe that 
$
Z = I_1 \oplus 0 
$
is the maximum rank slack; indeed  i) $Z = B - A_1 - A_2 - A_3,$ so it is a slack, and ii) 
the matrix 
\beq \label{UU}
U = \bpx 0 & 0 & 0 \\ 0 & 10 & 3 \\ 0 & 3 & 1 \epx
\eeq
satisfies $B \bullet U = A_i \bullet U = 0$ for all $i.$ Hence $U$ is 
 orthogonal to any slack matrix, so the rank of any slack matrix is at most $1.$ 

If $\alpha \neq 1, \,$ then \eref{lastsystem} is badly behaved; as proof, observe that   
$$
V \, := \, A_3 - A_2 - A_1 \, = \, \bpx 1 & 0 & \alpha - 1  \\ 
0 & 0 & 0 \\
\alpha - 1 & 0 & 0 \epx 
$$ 
is a certificate matrix as required by Theorem \ref{badsdp}.  

If $\alpha = 1, \,$ then \eref{lastsystem} is well behaved, and we can verify this using Theorem \ref{goodsdp} as follows. The $U$ matrix in \eref{UU} satisfies condition \eref{cond1}
of Theorem \ref{goodsdp}. As to condition \eref{cond2},  if the lower 
right $2 \ti 2$ block of $V := \sum_{i=1}^3 \lambda_i A_i + \mu B$ is zero, then $(\lambda_1, \lambda_2, \lambda_3, \mu)$ must be a linear combination of 
$$
(0,0,2,-1) \, {\rm and} \, (5,5,-1,-2),
$$
so for all such $(\lambda_1, \lambda_2, \lambda_3, \mu)$ the upper left $1 \times 2$ block of $V$ is also zero. 
 }
\eex
We return to Examples \ref{ex1}--\ref{alpha-ex} in Section \ref{npcapconp}.
As we will see there, the bad or good behavior of semidefinite systems 
can be verified using only an  elementary 
linear algebraic argument, 
{\em without ever referring to Theorems \ref{badsdp} or \ref{goodsdp}.} 
We will use Examples \ref{ex1}--\ref{alpha-ex} as illustrations. 

The reader may find it interesting to spot the $Z$ and $V$ excluded matrices in
 other pathological SDPs in the literature, e.g., in the 
instances in \cite{BentalNem:01, BonnShap:00, Barvinok:2002, VanBo:96, Todd:00, Ramana:97, TunWolk:12, LuoSturmZhang:97}.

Theorems \ref{badsdp} and 
\ref{goodsdp} simply follow from Theorem \ref{unif-d-thm} and from 
Lemma \ref{psd-dir} below, which describes 
the set of feasible directions and related sets in the 
semidefinite cone:  
\ble \label{psd-dir}
Let $Z$ be as in Assumption \ref{slack-ass}, and  
recall the definition of the set of feasible directions, and related sets 
from \eref{dirldir}-\eref{deftan}.
Then 
\beqa 
\label{psdldir}
\ldir(Z, \psd{n}) & = & \sym{r} \oplus \{ 0 \}, \\
\label{psdcldir}
\cl \dir(Z, \psd{n}) & = & \left\{ \,  \bpx Y_{11} & Y_{12} \\ Y_{12}^T & Y_{22} \epx \, \bigr| \bigl. \, Y_{22} \in \psd{n-r} \; \right\}, \\ 
\label{psdtan}
\tan(Z, \psd{n}) & = & \left\{ \, \bpx Y_{11} &  Y_{12} \\ Y_{12}^T & 0  \epx \, \bigr| \bigl. \, Y_{11} \, \in \, \sym{r} \, \right\}, \\ 
\label{psdcldirdir}
\dir(Z, \psd{n}) & = & \, \left\{ \, \bpx Y_{11} & Y_{12} \\ Y_{12}^T & Y_{22} \epx \, \bigr| \bigl. \, Y_{22} \in \psd{n-r}, \R(Y_{12}^T) \subseteq \R(Y_{22}) \; \right\}.
\eeqa
\ele
\qed

The proof of Lemma \ref{psd-dir} is given in Appendix B.

\pf{of Theorem \ref{badsdp}} 
By condition (1) of Theorem \ref{unif-d-thm} we see that 
\eref{p-sd} is badly behaved, iff 
there is a matrix $V \in \lin \{ A_1, \dots, A_m, B \}  \,$ 
such that 
\beq \nonumber \label{Vpsd} 
V \in \cl \dir(Z, \psd{n}) \setminus \dir(Z, \psd{n}).
\eeq
Thus our result follows from parts \eref{psdcldir} and 
\eref{psdcldirdir} in Lemma \ref{psd-dir}.
\qed

\pf{of Theorem \ref{goodsdp}} 
\co{To prove our result, we invoke the 
second part of  Theorem \ref{unif-d-thm}. }
We apply Theorem \ref{unif-d-thm} to the system (\ref{p-sd}). 
We first observe that a matrix $U \succeq 0$ is 
strictly complementary to $Z$ if and only if 
\beqa \nonumber 
U & = & \bpx 0 & 0 \\ 0 & U_{22} \epx, \; \mathrm{with} \; U_{22} \in \sym{n-r}, \, U_{22} \succ 0.
\eeqa
Next we note that 
 the first part of condition $\eref{tang}$ in Theorem \ref{unif-d-thm} 
holds iff 
there is such a $U$ that satisfies \eref{AiU}. 
By  \eref{psdldir} and \eref{psdtan} in Lemma \ref{psd-dir} 
the second part of condition \eref{tang} in Theorem \ref{unif-d-thm} 
holds iff all 
$V \in \lin \{ A_1, \dots, A_m, B \}$ which are of the form 
\beq \nonumber
V \, = \, \bpx V_{11}   & V_{12} \\
           V_{12}^T & 0  \epx
\eeq
satisfy $V_{12} = 0.$ This completes the proof.
\qed

To summarize, Theorems \ref{badsdp} and \ref{goodsdp} are a "combinatorial 
version" of Theorem \ref{unif-d-thm}.

We  note that  for semidefinite systems that are strictly feasible, 
a matrix similar to the $V$ matrix in Theorem \ref{badsdp} 
can make sure that the optimal primal-dual solution pair 
fails strict complementarity; see \cite{WeiWolk:06}. 

Although we focus on feasible systems, we obtain  
natural corollaries  about 
{\em weakly infeasible} SDPs, a class of pathological infeasible SDPs.
To describe the connection, note that the alternative system 
\beq \label{strinf}
Y \succeq 0, \, A_i \bullet Y = 0 \, (i=1, \dots, m), \, B \bullet Y = -1 
\eeq
gives a natural proof of infeasibility of \eref{p-sd}: if \eref{strinf} 
is feasible, then \eref{p-sd} is trivially infeasible.
However, \eref{p-sd} and \eref{strinf} may both be infeasible, in which case we call the semidefinite system
\eref{p-sd} {\em weakly infeasible.}

As background on weakly infeasible SDPs, we mention
 that Waki \cite{Waki:12} recently described a method for generating weakly infeasible SDPs based on
Lasserre's relaxation for polynomial optimization problems; Klep and Schweighofer \cite{KlepSchw:13} analyzed weakly infeasible SDPs
using real algebraic geometry techniques; and 
Lourenco et al \cite{Lourenco:13} proved that any weakly infeasible SDP with order $n$ matrices has 
a weakly infeasible subsystem with dimension at most $n-1.$   

To apply our machinery to weakly infeasible SDPs, we homogenize (\ref{p-sd}) to obtain the system
\beq \label{homogenized} 
\sum_{i=1}^m x_i A_i - x_0 B \preceq 0.
\eeq
\co{Let us assume that the system (\ref{homogenized})  satisfies Assumption 
\ref{slack-ass}.
First, suppose that \eref{p-sd} is weakly infeasible.
Then (\ref{homogenized}) is badly behaved, since  the objective function 
$\sup - x_0$ gives  a $0$ optimal value over this system, 
but there is no feasible dual solution 
(such a dual solution would be feasible in \eref{strinf}). 
Hence by Theorem \ref{badsdp} the excluded matrices $Z$ and $V$ appear 
in (\ref{homogenized}).    In turn, if \eref{homogenized} satisfies the conditions of
Theorem \ref{goodsdp} and hence  it is well behaved,
then (\ref{p-sd})   cannot be weakly infeasible. }
Assume that the system (\ref{homogenized})  satisfies Assumption 
	\ref{slack-ass}.
	First, suppose that \eref{p-sd} is weakly infeasible.
	Then (\ref{homogenized}) is badly behaved, since 
	\beq \label{ffff}
	\sup \{ \, x_0 \, | \, (x,x_0) \, {\rm is \, feasible \, in \, } (\ref{homogenized}) \} = 0,
	\eeq
	but there is no solution feasible in the dual of (\ref{ffff}) 
	(such a dual solution would be feasible in \eref{strinf}). 
	Hence by Theorem \ref{badsdp} the excluded matrices $Z$ and $V$ appear 
	in (\ref{homogenized}).    In turn, if \eref{homogenized} satisfies the conditions of
	Theorem \ref{goodsdp} and hence  it is well behaved,
	\co{hence the characterization of Theorem 
		(\ref{goodsdp}) applies, }
	then (\ref{p-sd})   cannot be weakly infeasible.

\section{Reformulations. Badly behaved semidefinite systems are in \mbox{$\np \cap \conp$}} 
\label{npcapconp}  

\subsection{Reformulations}

To motivate the discussion of this section, we recall a basic result from the theory of linear equations:
\bit
\item[] ``The system $Ax=b$ is infeasible if and only if its row echelon form contains the 
equation  $ \la 0, x \ra = \alpha,$ where $\alpha \neq 0.$''
\eit
Since the "if" direction is trivial, we  will --- informally --- say that the 
 row echelon form is an {\em easy-to-verify certificate, or witness},  of infeasibility.

In this section we describe  analogous results for a very different problem :
we show how to transform \eref{p-sd} into an equivalent system whose bad or good behavior is 
trivial to verify. 
As a  corollary we  prove that  badly (and  well) behaved semidefinite systems are in 
$\np \cap \conp$ in the  
real number model of computing. (In this model we can store arbitrary real numbers in unit space 
and perform arithmetic operations in unit time; 
see e.g. \cite{Cuckeretal:98}. We do not claim that badly behaved semidefinite systems are in ${\cal P},$ i.e., we do not 
provide a polynomial time algorithm to decide whether \eref{p-sd} is badly behaved. We discuss this point in more detail at the end of this section.) 

\co{
Observe that Theorem \ref{badsdp} already provides the excluded matrices $Z$ and $V$ as witnesses 
 of bad behavior.
However, $Z$ and $V$ are not easy-to-verify certificates: to convince a reader that they lead to bad 
behavior, we must refer to the proof of the "if" direction of 
Theorem \ref{badsdp}, which is lengthy.
These matrices, by themselves, are also not polynomial time verifiable certificates: 
a verifier  would also need an accompanying  witness
 that $Z$ is indeed a maximum rank slack (say if it has rank 5, then we need a proof that there is no rank 6 feasible slack). 
}

\co{
Note that Theorem \ref{badsdp} already provides the excluded matrices $Z$ and $V$ as witnesses 
of bad behavior.
However, these witnesses are not easy-to-verify:  to convince a reader that they lead to bad 
behavior, we must invoke the proof of the "if" direction of 
Theorem \ref{badsdp}, which is lengthy.
These matrices, by themselves, also do not give a polynomial time verifiable proof 
that (\ref{p-sd}) is badly behaved: such  a proof would  need 
an accompanying  witness that $Z$ is indeed a maximum rank slack (say if it has rank 5, then we must prove  that there is no rank 6 feasible slack). 
}

\co{
Observe that Theorem \ref{badsdp} already provides the excluded matrices $Z$ and $V$ as witnesses 
of bad behavior.
However, these witnesses are not easy-to-verify:  to convince a reader that they lead to bad 
behavior, we must invoke the proof of the "if" direction of 
Theorem \ref{badsdp}, which is lengthy.
These matrices, by themselves, also do not give a polynomial time verifiable proof 
that (\ref{p-sd}) is badly behaved. The reason is that -- until now -- 
we have not shown a way to prove  such a proof would  need 
an accompanying  certificate that $Z$ is indeed a maximum rank slack (say if it has rank 5, then we must prove  that there is no rank 6 feasible slack). 
}

We first  define the type of transformation that we use on \eref{p-sd}.

\bdef \label{reform-def} {\rm We obtain an {\em elementary reformulation,} or simply a 
	{\em  reformulation}, 
 of $(\sdpc)$ 
by a sequence of the following operations:
\benum
\item \label{rotate} Apply a 
rotation $T^T()T$ to all $A_i$ and $B, \,$ where $T = I_r \oplus M \,$ and 
$M$ is invertible.
\item \label{slack} Replace $B$ by $B + \sum_{j=1}^m \mu_j A_j, \,$ where $\mu \in \rad{m}.$ 
\item \label{exch} Exchange $(A_i, c_i)$ and $(A_j, c_j), \,$ where $ i \neq j.$ 
\item \label{trans} Replace $(A_i, c_i)$ 
by $(\sum_{j=1}^m \lambda_j A_j, \sum_{j=1}^m \lambda_j c_j),  \,$ where 
$ \lambda \in \rad{m}, \, \lambda_i \neq 0.$ 
\eenum
We obtain an elementary reformulation of the system \eref{p-sd} by applying the preceding 
operations 
with some $c.$ 
}
\end{Definition}
Clearly, in all reformulations of \eref{p-sd} the maximum rank slack 
is the same. 

Where do these operations come from?
Operations \eref{exch} and \eref{trans}  are equivalent to
elementary row operations (inherited 
from Gaussian elimination) done on $(\sddc):$ 
\bit
\item Operation \eref{exch}  exchanges  the dual equations
$
A_i \bullet Y = c_i \, {\rm and} \, A_j \bullet Y = c_j;  
$ and 
\item Operation \eref{trans} replaces
the dual equation
$
A_i \bullet Y = c_i \, {\rm by} \, \sum_{j=1}^m (\lambda_j A_j) \bullet Y = \sum_{j=1}^m \lambda_j c_j.  
$
\eit
\ble \label{ref-lemma}
The system \eref{p-sd} is well behaved if and only if its elementary reformulations are.
\ele

\pf{} Operations (1)-(4) of Definition \ref{reform-def} 
keep the value of $(\sdpc)$ finite, if it is finite; and infinite, if it is infinite. 
Suppose now that $Y$ is feasible in $(\sddc)$ with value, say,  $\alpha, $ and 
we apply operations (1) and (2) with rotation matrix $T \,$ and vector $\mu.$ 
Then identity \eref{txt} implies that  
$T^{-T}YT$ is feasible in the dual of the reformulated problem with value 
$\alpha + \sum_{j=1}^m \mu_j c_j.$ Operations (3) and (4) preserve the 
feasibility and objective value of a solution of $(\sddc).$
Thus if \eref{p-sd} is well behaved, so are its 
reformulations, and this completes the proof of the "Only if" direction.
Since (\ref{p-sd}) is  a reformulation of its reformulations, the "If" direction follows as well. 
\qed

\subsection{Reformulating (\ref{p-sd}) to verify maximality of the maximum rank slack}

Recall that $Z$ is the maximum rank slack in \eref{p-sd} described  
 in Assumption \ref{slack-ass}. 
We reformulate \eref{p-sd} in two steps. In the first step, given 
in Lemma \ref{max-sl-lemma}, we reformulate \eref{p-sd} so the resulting system has 
easy-to-verify witnesses  that $Z$ is a maximum rank slack.
(The $Y_j$ matrices in Lemma \ref{max-sl-lemma} will be the witnesses.)

In Lemma \ref{max-sl-lemma} we rely on a facial reduction algorithm (see 
\cite{BorWolk:81, BorWolk:81B, WakiMura:12,Pataki:13}).
It is important that in Lemma \ref{max-sl-lemma} we only use 
rotations, i.e., type \eref{rotate} operations
of  Definition \ref{reform-def}. 

\ble \label{max-sl-lemma}
The system \eref{p-sd} 
has a reformulation 
\beq \label{p-sd-p} \tag{\mbox{$P_{SD}'$}}
\sum_{i=1}^m x_i A_i' \preceq B' 
\eeq
and there exist symmetric matrices of the form
\co{
	$$
Y_{j} = 
\bordermatrix{& \hspace{1cm}  &   \overbrace{\qquad}^{\textstyle{r_j}} &  \overbrace{\qquad}^{\textstyle{r_{j-1} + \dots + r_1}} \cr
	& 0   &  0   &  \ti \cr
	& 0   &  I   &  \ti \cr
	& \ti   &  \ti  &  \ti \cr} \, (i=1, \dots, \ell) 
$$
}
\beq \label{Yjform} 
Y_{j} = 
\bordermatrix{& \hspace{1cm}  &   \overbrace{\qquad}^{\scriptstyle{r_j}} &  \overbrace{\qquad}^{\scriptstyle{r_{j-1} + \dots + r_1}} \cr
	& 0   &  0   &  \ti \cr
	& 0   &  I   &  \ti \cr
	& \ti   &  \ti  &  \ti \cr} \, (j=1, \dots, \ell) 
\eeq
where $\ell \geq 0, r_1 > 0, \dots, r_\ell > 0, \, $ 
$r_1 + \dots + r_\ell = n - r, \,$  and  
\beq \label{YjS} 
Y_j \bullet B' = Y_j \bullet A_i' = 0 
\eeq
holds for all $i$ and $j.$ 
Here the 
$\ti$ symbols denote blocks with arbitrary elements in  the $Y_j$ matrices. 
\ele
\qed

If $Z = I, \,$ i.e., (\ref{p-sd}) satisfies Slater's condition, then we just take 
$B' = B, \, A_i' = A_i$ for all $i$ and 
$\ell = 0$ in Lemma \ref{max-sl-lemma}. 

To build intuition, we first establish  why the $Y_j$ matrices  
indeed prove that the rank of any slack matrix is at most $r.$
Let  $S$ be a slack in (\ref{p-sd}), and $Y_1, \dots, Y_\ell$ 
as in the statement of
Lemma \ref{max-sl-lemma}.  Then $S = B' - \sum_i x_i A_i' \,$ for some $x \in \rad{m}.$ 
So 
$Y_1 \bullet S = 0$ and $S \succeq 0, \,$ hence  the last $r_1$ rows and columns 
of $S$ are zero;  $Y_2 \bullet S = 0$ and $S \succeq 0$ imply that 
the next $r_2$ rows and columns of $S$ are zero, and so on. Inductively we find that the last $r_1 + \dots + r_\ell = n-r$ rows and columns of $S$ are zero, hence $S$ must have rank at most $r.$  

Thus we can prove that $Z$ is a maximum rank slack 
 in \eref{p-sd-p} (hence also in (\ref{p-sd}))   using
\benum
\item a vector $x \in \rad{m}$ such that $Z = B' - \sum_{i=1}^m x_i A_i', $ and 
\item the $Y_j$ matrices of Lemma \ref{max-sl-lemma}.
\eenum 
We next illustrate Lemma \ref{max-sl-lemma}.

\bex{\rm (Examples  \ref{ex1}, \ref{ex2}, \ref{ex3} and \ref{alpha-ex}  continued)
	In all these examples it is easy to show why the maximum rank slack is indeed a slack. Also, in Example \ref{ex1} 
$$ 
Y_1 = \bpx 0 & 0 \\ 0 & 1 \epx
$$
is orthogonal to all constraint matrices (using  the $\bullet$ inner product), 
so it proves that the rank of any slack matrix is at most one. 

In Example \ref{ex2} the matrix $Y_1 = 0 \oplus I_1$ proves that the rank of any slack is at most one, and 
in Example  \ref{ex3}  the matrix 
$Y_1 = 0 \oplus I_2 $ proves that the rank of 
any slack is at most two.  
(So the first three examples do not even need to be reformulated to have a convenient proof that $Z$ is a maximum rank slack.)

In Example \ref{alpha-ex} we  let
$$
T = \bpx 1 & 0 & 0 \\ 0 & 1 & 3 \\ 0 & 0 & 1 \epx,
$$
and apply the rotation $T^T()T$ to all matrices to obtain the system 
\beq \label{alpha-cont} 
x_1 \bpx 0   &  0 &    1 \\
0 &    1  &   0 \\
1  &   0 &   -1
\epx + x_2 \bpx 0   &  1 &    0 \\
1 &    0 &    1 \\
0  &   1 &    0 \epx +
x_3 \bpx 1 & 1 & \alpha \\ 1 & 1 & 1  \\ \alpha & 1 & -1 \epx \preceq  \bpx 2 & 2 & \alpha +1 \\ 2 & 2 & 2 \\ \alpha + 1 & 2 &  -2 \epx.
\eeq
 Now  $Y_1 = 0 \oplus I_2$ is orthogonal to all constraint matrices in 
 (\ref{alpha-cont}), and this proves that the rank of any slack is at most one, so 
 $Z = I_1 \oplus 0$ is a maximum rank slack. 
 
 In Appendix A we give a larger example, in which we need two $Y_j$ matrices to prove
 that any slack matrix has rank at most $2.$ 
}
\eex
\pf{of Lemma \ref{max-sl-lemma}} 
\co{We first note that to satisfy equation (\ref{YjS}) it is enough to have 
$Y_j \bullet A_i' = Y_j \bullet B = 0 \,$ for all $i$ and $j.$ Indeed, if $S$ is any slack matrix, then 
$Y_j \bullet S = 0$ follows, since $S = B' - \sum_i x_i A_i'.$ }

To find the reformulation assume that $k \geq 0,$ 
we have a reformulation of the form (\ref{p-sd-p}) and matrices $Y_1, \dots, Y_k$ such that 
\eref{YjS} holds for all $i$ and for $j=1, \dots, k.$ 
At the start $k=0$ and $B' = B, \, A_i' = A_i \,$ for all $i.$ 
For brevity, let $s_k = r_1 + \dots + r_k. $  We claim 
that 
$$
s_k \leq n - r \, {\rm holds.}
$$ 
This indeed follows since if $S \succeq 0$ is a slack in 
(\ref{p-sd-p}) then  (using the same argument that we used before) 
the last $s_k$ rows and columns of $S$ must be zero. 

If $s_k = n - r,$  
we set $\ell = k, \,$ and stop; otherwise, we define the cone $K = \psd{n} \cap Y_1^\perp \dots \cap Y_k^\perp.$ 
Clearly, $K$ and its dual cone $K^*$ are of the form 
$$
K \, = \, \left\{ \,  \bpx Y_{11} & 0  \\ 0  & 0 \epx \, \bigr| \bigl. \, Y_{11} \in \psd{n - s_k} \; \right\}, 
\, K^* \, = \, \left\{ \,  \bpx Y_{11} & Y_{12} \\ Y_{12}^T & Y_{22} \epx \, \bigr| \bigl. \, Y_{11} \in \psd{n - s_k} \; \right\}.
$$
Next, define the affine subspace 
$$
H \, = \, \lin \, \{ \, A_1', \dots, A_m' \, \} + B'.
$$
Since $Z$ is also a maximum rank slack in 
$\psdsysp$, and $r < n - s_k, \,$ we have $H \cap K \neq \emptyset, \,$ 
$H \cap \ri K = \emptyset, \,$ hence 
\mbox{$H^\perp \cap (K^* \setminus K^\perp) \neq \emptyset$} 
by a classic theorem of the alternative 
(see e.g. Lemma 1 in \cite{Pataki:13}).

Let 
$$
Y_{k+1} \in H^\perp \cap (K^* \setminus K^\perp).
$$
Since $Y_{k+1} \bullet Z = 0, \,$ we have 
\co{
	$$
Y_{k+1} = 
\bordermatrix{& \scr{r}   &  \hspace{1cm} &  \scr{s_i} \cr
& 0   &  0   &  \ti \cr
& 0   &  Y'   &  \ti \cr
& \ti   &  \ti  &  \ti \cr}
$$
}
$$
Y_{k+1} = 
\bordermatrix{& \overbrace{\qquad}^{\scriptstyle{r}}   &  \hspace{1cm} &  \overbrace{\qquad}^{\scriptstyle{s_{k}}}  \cr
	& 0   &  0   &  \ti \cr
	& 0   &  Y'   &  \ti \cr
	& \ti   &  \ti  &  \ti \cr}
$$
for some $Y' \succeq 0.$ (Again, the $\ti$ symbols stand for submatrices with arbitrary elements). 
Let $r_{k+1}$ be the number of positive eigenvalues of $Y';$ since 
$Y_{k+1} \not\in K^\perp, \, $ we have $r_{k+1} > 0.$ 

Let $Q$ be an invertible matrix such that 
$Q^T Y' Q = 0 \oplus I_{r_{k+1}},$ 
and $T = I_r \oplus Q \oplus I_{s_k}.$ 
We apply the rotation $T^T()T$ to $Y_1, \dots, Y_{k+1}, \,$ and the rotation 
$T^{-1}()T^{-T}$ to all $A_i'$ and to $B'.$  

By \eref{txt} the equation (\ref{YjS}) holds for all $i$ and for 
$j=1, \dots, k+1.$  
By the form of $T$ now  $Y_1, \dots, Y_{k+1}$ are in the required shape 
(see equation (\ref{Yjform})).  
We then set $k = k+1$ and 
continue. 

Clearly, our algorithm terminates in finitely many steps, so the proof is complete.
\qed

\subsection{Reformulating (\ref{p-sd}) to verify that it is badly behaved}

 In Theorem \ref{badsdp-re} we give the {\em final} reformulation of
\eref{p-sd} to prove its bad behavior. 
We  point out that in  Theorem \ref{badsdp-re} the proof of the "if" direction is  elementary, thus 
the reformulated system \eref{p-sd-bad} 
is  an easy-to-verify certificate that \eref{p-sd} is badly behaved. 

\bth \label{badsdp-re}
The system \eref{p-sd} is badly behaved if and only if it has a reformulation 
\beq \label{p-sd-bad} \tag{$P_{\text{SD,bad}}$}
\sum_{i=1}^k x_i \bpx F_i & 0 \\ 0 & 0 \epx + \sum_{i=k+1}^m x_i \bpx F_i & G_i \\ G_i^T & H_i \epx \, \preceq \, \bpx I_r & 0 \\ 0 & 0 \epx = Z,
\eeq
where 
\benum
\item matrix $Z$ is the maximum rank slack, and its maximality can be  verified by matrices 
$Y_1, \dots, Y_\ell, \,$ as given by Lemma \ref{max-sl-lemma}.
\item \label{indep} The matrices 
$$
\bpx G_i \\ H_i \epx \, (i=k+1, \dots, m) 
$$
are linearly independent.
\item $H_m \succeq 0.$
\eenum
\enth 
\pf{(If)} By Lemma \ref{ref-lemma} 
it is enough to prove that \eref{p-sd-bad} is badly behaved. Let  
$x$ be feasible in $\psdsysb$  with a corresponding slack $S.$ 
Note that 
the last $n-r$ rows and columns 
of $S$ must be zero, otherwise $\frac{1}{2}(S+Z)$ would be
a slack with larger rank than $Z.$ 
Hence, by condition (2) we must have 
$x_{k+1} \, = \, \dots \, = \, x_m \, = \, 0.$ 
Next, let us consider the SDP
\beq \label{mofo}
\sup \, \{ \, - x_m \, | \, x \, \text{is feasible in } \mbox{\eref{p-sd-bad}} \, \},
\eeq
which, by the above argument,  has optimal value $0.$ 
We prove that its dual cannot have a feasible solution with value $0, \,$ 
so suppose that 
$$
Y \, = \, \bpx Y_{11} & Y_{12} \\ Y_{12}^T & Y_{22} \epx \succeq 0
$$
is such a solution. By $Y \bullet Z = 0$ we get $Y_{11} = 0, \,$ hence by psdness of 
$Y$ we deduce $Y_{12} = 0.$ 
Thus 
$$
\bpx F_m & G_m \\ G_m^T & H_m \epx \bullet Y = H_m \bullet Y_{22} \geq 0,
$$
which contradicts the assumption that $Y$ is feasible in the dual of (\ref{mofo}). 

\pf{(Only if)} 
We start with the system (\ref{p-sd-p}) given by Lemma 
\ref{max-sl-lemma} and further reformulate it. 
For brevity we denote the constraint matrices on the left hand side by $A_i'$ throughout the process.

We first replace $B'$ by $Z$ in (\ref{p-sd-p}). Since the resulting system is still 
badly behaved, by Theorem \ref{badsdp} 
there is a matrix of the form 
$$
V \, = \, \lambda_0 Z + \sum_{i=1}^m \lambda_i A_i' \, = \, \bpx V_{11} & V_{12} \\
                V_{12}^T & V_{22} 
\epx, \, 
$$
with 
$V_{11} \in \sym{r}, V_{22} \succeq 0, \, {\rm and} \, \R(V_{12}^T) \not \subseteq \R(V_{22}).$
By the form of $Z$ we can assume $\lambda_0 = 0$ (otherwise we can replace $V$ by $V - \lambda_0 Z$). 

Note that the block of $V$ comprising the last $n-r$ columns must be nonzero.
We pick  an $i \,$ such that $\lambda_i \neq 0, \,$  
replace $A_i'$ by $V,$  
then switch $A_i'$ and $A_m'.$ 
Next we choose a maximal subset of the $A_i'$ matrices so their blocks
comprising the last $n-r$ columns are linearly independent. We let $A_m'$ to be one 
of these matrices (this can be done, since $A_m'$ is now the $V$ certificate matrix), and permute the $A_i'$ so this special subset 
becomes $A_{k+1}', \dots, A_m'$ for some $k \geq 0.$ 

We finally add suitable multiples of $A_{k+1}', \dots, A_m'$ to $A_1', \dots, A_k'$ 
to zero 
out the last $n-r$ columns and rows of the latter, 
and arrive at the required reformulation.
\qed

\bex{\rm (Examples  \ref{ex1}, \ref{ex2} and \ref{alpha-ex}  continued)
	The first two of these examples are already in the standard form
	(\ref{p-sd-bad}).
	Suppose now $\alpha \neq 1 \,$ in  Example  \ref{alpha-ex},
	i.e., the system \eref{lastsystem} is badly behaved.  	
	Recall that by a rotation we brought \eref{lastsystem} to the simpler form 
	\eref{alpha-cont}. Then in \eref{alpha-cont} we set
	\beqast 
	B & := & B - A_1 - A_2 - A_3, \\
	A_3 & := & A_3 - A_1 - A_2,
	\eeqast
	and obtain the system 
	\beq \label{lastsystemre}
	\ba{rl}
	& x_1 \bpx 0 & 0 & 1 \\ 0 & 1 & 0 \\ 1 & 0 & -1 \epx 
	+ x_2 \bpx 0 & 1 & 0 \\ 1 & 0 & 1 \\ 0 & 1 & 0  \epx
	+ x_3 \bpx 1 & 0 & \alpha - 1 \\ 0 & 0 & 0 \\ \alpha-1 & 0 & 0 \epx
	\preceq 
	\bpx 1 & 0 & 0 \\ 0 & 0 & 0 \\ 0 & 0 & 0 \epx, 
	\ena
	\eeq
	which is in the standard form  \eref{p-sd-bad} (with $k=0$). 
	The objective function $\sup - x_3$ 
   yields a zero optimal value over (\ref{lastsystemre}) 
	but there is no dual solution with the same value: we can argue this as in the 
	proof of the "if" direction in Theorem \ref{badsdp-re}. 
}
\eex

Note that the certificate matrix $V$ of Theorem \ref{badsdp} appears in the system
(\ref{p-sd-bad}) as the last matrix on the left hand side.

\subsection{Reformulating (\ref{p-sd}) to verify that it is well  behaved}

We now turn to well behaved semidefinite systems, and 
in Theorem \ref{goodsdp-re} we show 
how to reformulate them to easily verify their good behavior.
In Theorem \ref{goodsdp-re} we also show block-diagonality of dual 
optimal solutions.  Note that the proof of the "if" direction in  Theorem \ref{goodsdp-re} is easy, so the system \eref{p-sd-good} is an easy-to-verify certificate of good behavior. 

\bth \label{goodsdp-re}
The system \eref{p-sd} is well behaved if and only if it has a reformulation 
\beq \label{p-sd-good} \tag{$P_{\text{SD,good}}$}
\sum_{i=1}^k x_i \bpx F_i & 0 \\ 0 & 0 \epx + \sum_{i=k+1}^m x_i \bpx F_i & G_i \\ G_i^T & H_i \epx \preceq \bpx I_r & 0 \\ 0 & 0 \epx = Z,
\eeq
where 
\benum
\item the matrix $Z$ is the maximum rank slack.
\item \label{indep2} The matrices 
$
H_i \, (i=k+1, \dots, m) 
$
are linearly independent.
\item \label{HiI} $H_{k+1} \bullet I = \dots = H_m \bullet I = 0.$ 
\eenum
Also, if \eref{p-sd} is well behaved, and 
the value of $(\sdpc)$ is finite, then there is an optimal 
dual matrix in $\psd{r} \oplus \psd{n-r}.$ 
\enth
\pf{ (If and block-diagonality)} 
Let  $c$ be such that 
\beq \label{origgood}
v := \sup \, \{ \, \sum_{i=1}^m c_i x_i \, | \,  x \, \text{is feasible in } \psdsysg \, \}
\eeq
is finite. 
By the proof of Lemma \ref{ref-lemma} 
it suffices to prove that the dual of 
(\ref{origgood}) 
has a block-diagonal solution with value $v.$ 
An argument like in the proof of Theorem \ref{badsdp-re} proves that 
$
x_{k+1} = \dots = x_m = 0 
$
holds for any  $x$ feasible in \eref{origgood}, so 
\beq \label{reduced}
v \, = \, \sup \, \{ \, \sum_{i=1}^k c_i x_i \, | \,  \sum_{i=1}^k x_i F_i \preceq I_r \, \}.
\eeq
Since \eref{reduced} satisfies Slater's condition, 
there is $Y_{11}$ feasible in its dual with $Y_{11} \bullet I_r = v.$ 

As the $H_i$ are linearly independent, we can choose $Y_{22} \in \sym{n-r}$ 
(which is possibly not psd)
such that 
$$
Y := \bpx Y_{11} & 0 \\ 0 & Y_{22} \epx
$$
satisfies the equality constraints of the dual of \eref{origgood}.
We then add a positive multiple 
of the identity to $Y_{22}$ to make $Y$ psd. 
Taking condition (3) into account we can see that after this operation
$Y$ is feasible in the dual of \eref{origgood} 
and clearly 
$Y \bullet Z = v$ holds. The proof is now complete. 

\pf{(Only if)} We again start with the system (\ref{p-sd-p}) that Lemma \ref{max-sl-lemma} provides; now (\ref{p-sd-p}) is well behaved.
(We also note that the $U$ matrix of Theorem 
\ref{goodsdp} became the $Y_1 = 0 \oplus I_{n-r}$ 
matrix of Lemma \ref{max-sl-lemma}, after we rotated it.)
We first replace $B'$ by $Z.$ 
Next we choose a maximal subset of the $A_i'$ whose 
lower principal $(n-r) \ti (n-r)$ blocks are linearly independent. 
We permute the $A_i'$ if needed, to make this subset 
$A_{k+1}', \dots, A_m'$ for some $k \geq 0.$ 

To complete the process we add multiples of $A_{k+1}', \dots, A_m'$  to 
$A_1', \dots, A_k'$ to zero out the lower principal  
$(n-r) \ti (n-r)$ block of the latter. By Theorem \ref{goodsdp} 
the upper right $r \ti (n-r)$  block of 
$A_1', \dots, A_k'$ 
and the symmetric counterpart also become 
zero. 
This concludes the proof.
\qed

\bex{\rm \label{continued} (Examples \ref{ex3} and  \ref{alpha-ex}  continued) 
	In Example \ref{ex3} 
	the system (\ref{seinfeld}) is already in the form of (\ref{p-sd-good}).
	
	Suppose now $\alpha =1  \,$ in Example \ref{alpha-ex}, 
	i.e., \eref{lastsystem} is well behaved.
Recall that we transformed this system into  the system 
\eref{lastsystemre} (in Example \ref{continued}; note that this can be done 
independently of the value of  $\alpha$). 
We then switch the first and third matrices in \eref{lastsystemre} to get 
	\beq \label{lastsystemrere}
	\ba{rl}
	& x_1 \bpx 1 & 0 & 0 \\ 0 & 0 & 0 \\ 0 & 0 & 0 \epx
	+ x_2 \bpx 0 & 1 & 0 \\ 1 & 0 & 1 \\ 0 & 1 & 0  \epx
	+ x_3 \bpx 0 & 0 & 1 \\ 0 & 1 & 0 \\ 1 & 0 & -1 \epx 
	\preceq 
	\bpx 1 & 0 & 0 \\ 0 & 0 & 0 \\ 0 & 0 & 0 \epx
	\ena
	\eeq
 in the standard form \eref{p-sd-good} (with $k=1$). 
}
\eex

We next discuss some implications of Theorem \ref{goodsdp-re}.
First, as the proof of the "if" direction shows, we 
can compute an optimal solution of \eref{origgood} from an optimal solution of 
the reduced problem \eref{reduced}; to do so, we only need  to solve 
a linear system 
of equations (to find $Y_{22}$) and do a linesearch (to make $Y_{22}$ psd).  

Second,  loosely speaking, the system  (\ref{p-sd-good}) can be partitioned into 
a strictly feasible part,  and a linear part, which  corresponds to variables $x_{k+1}, \dots, x_m.$ 

Third, how do we generate a well behaved semidefinite system? 
Theorem \ref{goodsdp-re} can help us   to do this: we can choose 
matrices $Z, H_i, G_i, F_i$ to obtain a system in the form (\ref{p-sd-good}), 
then arbitrarily reformulate it, while keeping it well behaved. 
In fact, according to Theorem \ref{goodsdp-re}, we can obtain {\em any} well behaved semidefinite system in this manner.

In related work, Bomze et al in \cite{BomzeSchach:12} describe methods 
to generate pathological conic LP instances  from other pathological conic LPs.
\co{to use conic programs
(in particular, SDPs and programs over the copositive cone) 
with certain pathological features, to generate other pathological conic programs. }
Their results  differ from ours, since they need to start with a pathological conic LP.  

We also note that using Lemma \ref{pataki-cl} the authors 
in Theorem  3.2 in \cite{Dima:15} 
characterized the situation when the projection of $\psd{n}$ 
onto some entries is closed; we can view Theorem \ref{goodsdp-re} as a generalization of this result.

\subsection{Badly behaved semidefinite systems are in $\np \cap \conp.$ Certificates to verify (non)closedness of the linear image of the semidefinite cone}

We now state our main complexity result: 
\bth \label{npconp} 
Badly (and well) behaved semidefinite systems are in 
$\np \cap \conp$ in the real number model of computing.
\enth
\pf{} We give the following certificates to check the status of \eref{p-sd}: 
(1) a reformulation of (\ref{p-sd}) into the form \eref{p-sd-bad} or \eref{p-sd-good}; (2) 
the $Y_j$ matrices of Lemma \ref{max-sl-lemma}
 to verify that $Z \,$ is indeed a maximum rank slack; 
(3) a matrix $T = I_r \oplus M, \,$ 
and $\mu \in \rad{m}, \,$ which were 
used to transform (\ref{p-sd}) into  \eref{p-sd-bad} or \eref{p-sd-good}. 

The verifier first checks that 
\eref{p-sd-bad} or \eref{p-sd-good} is indeed a reformulation of \eref{p-sd}; 
then verifies the properties of \eref{p-sd-bad} or \eref{p-sd-good} as given in Theorems 
\ref{badsdp-re} or \ref{goodsdp-re};  then the proof 
of the ``If'' part in Theorems \ref{badsdp-re} or \ref{goodsdp-re} shows that these systems are well- or badly behaved.
\qed

Assume that we are working with the real number model of computing. We don't claim to have a polynomial time algorithm to decide whether (\ref{p-sd}) 
is badly behaved; in particular, we don't have a polynomial time algorithm to compute the $Z$ and $V$ excluded matrices of Theorem \ref{badsdp},
  or one to compute the reformulated systems 
(\ref{p-sd-bad}) or (\ref{p-sd-good}). 

In analogy, if (\ref{p-sd}) is feasible, we can verify this in polynomial time (by plugging in a feasible $x$).
If (\ref{p-sd}) is infeasible, we can also verify this in polynomial time, using one of the infeasibility certificates in 
\cite{Ramana:97, KlepSchw:13, WakiMura:12, LiuPataki:15}.
However, we don't know how to decide in polynomial time whether (\ref{p-sd}) is feasible.

	 	Thus feasibility of a semidefinite system is similar to the bad behavior of a feasible system: both properties are in 
	 	$\np \cap \conp, \,$ but neither is known to be in ${\cal P}.$

To conclude this section, we briefly discuss easy-to-verify certificates for the (non)closedness 
of the linear image of $\psd{n}.$ 
All linear maps that map from $\sym{n}$ to $\rad{m}$ are of the form  
 $\A^*: \sym{n} \rightarrow \rad{m}, \,$ 
where 
$$
\A(x) = \sum_{i=1}^m x_i A_i, \, \A^*(Y) = (A_1 \bullet Y, \dots, A_m \bullet Y)^T
$$
and $A_i \in \sym{n}$ for all $i.$ 
We know that  $\A^*(\psd{n})$ is closed if and only if the homogeneous 
system 
\beq \label{homogenized2} 
\sum_{i=1}^m x_i A_i \preceq 0
\eeq
is well behaved (this is immediate from Lemma \ref{duffin}). 
\co{
	Thus we obtain easy-to-verify certificates for the (non)closedness of
	$\A^*(\psd{n})$ if we reformulate the system (\ref{homogenized2}) 
	into the standard form of (\ref{p-sd-bad}) or (\ref{p-sd-good}). .}
Thus  reformulating this homogeneous system 
into the standard forms
of \eref{p-sd-bad} or  \eref{p-sd-good}  gives 
easy-to-verify certificates of the closedness or nonclosedness of $\A^*(\psd{n}).$ 

	To illustrate this point we revisit Examples \ref{1stex} and  \ref{MCex}. 
	The semidefinite system
	\beq \label{MCeq} 
	- \M(x) \preceq 0,
	\eeq 
	where $\M$ is the linear map defined there, is badly behaved 
	(since the image of the semidefinite cone under $\M^*$ is not closed). We can apply 
	the machinery of this paper to study the system (\ref{MCeq}); 
	e.g., we can find the $Z$ and $V$ excluded matrices of Theorem 
	\ref{badsdp}, and reformulate (\ref{MCeq}) into the standard form (\ref{p-sd-bad}). We leave the details to the reader.

\section{Concluding remarks}
\label{section-conclude}

Theorem \ref{badsdp} gives the $Z$ and $V$ excluded matrices 
to characterize bad behavior of (\ref{p-sd}). We can carry this idea further, and 
prove  the following result: 
\bcor \label{badsdp-cor}
Suppose that in addition to the operations of Definition \ref{reform-def} we allow 
a sequence of the following operations: 
\benum
\item \label{delrow} 
Delete row $i$ and column $i$ from all matrices, where $i \in \{1, \dots, n \}.$
\item Delete a constraint matrix.
\eenum
Then we can bring any badly  behaved semidefinite system to the form of 
\beq \label{ex1-alpha-2}
\ba{rl}
x_1 \bpx \alpha & 1 \\ 1 & 0 \epx \preceq \bpx 1 & 0 \\ 0 & 0 \epx,
\ena
\eeq
where $\alpha$ is some real number.
\ecor
\pf{} Suppose that \eref{p-sd} is badly behaved and let us recall the form of the maximum rank slack in 
Assumption \ref{slack-ass}. 
We first add multiples of the $A_i$ to $B$ to make sure 
that the right hand side is the maximum rank slack. 
Next we let $V$ to be a certificate matrix as given by Theorem \ref{badsdp}; we can assume that $V$  
is the linear combination of 
the $A_i$ only; we reformulate, so $V$ becomes a constraint matrix. 

As we show in Lemma \ref{psd-dir}, we can apply a rotation $T^T()T$ to $V$ (where $T = I_r \oplus M \,$ for some invertible $M$)  
to bring $V$  to the form
\beq \label{newV}
V \, = \, \bpx V_{11}   & V_{12} & V_{13} \\
               V_{12}^T & I_s    & 0      \\
               V_{13}^T &  0     & 0 
                 \epx,
\eeq
where $V_{11}$ is $r \ti r, \,$ $s \geq 0$ 
and $V_{13} \neq 0.$ We apply the rotation $T^T()T$  to all constraint matrices, and after this operation
 $V$ is of 
the form specified in (\ref{newV}). 
Suppose now that $v_{ij} \neq 0, \,$ where $1 \leq i \leq r$ and 
$ r + s + 1 \leq j \leq n.$ We rescale $V$ to make sure that  $v_{ij} = 1 \,$ holds, then 
delete all rows and columns from the 
constraint matrices whose index is not $i$ nor $j, \,$ to obtain system \eref{ex1-alpha-2}.  
\qed

Excluded minor results in graph theory, such as Kuratowski's theorem,
show that a graph lacks a certain fundamental property, if and only if 
it can be reduced to a minimal such 
graph by a sequence of elementary operations.
Corollary \ref{badsdp-cor} resembles such results, since system \eref{ex1-alpha-2} is 
trivially badly behaved.

We can define the well- or badly behaved nature of conic linear systems 
in a different form,  and 
characterize such systems. For instance, we call the dual system 
\beq \label{dform}
\A^* y = c, \, y \in K^*,
\eeq
well behaved, if for all $b$ 
dual objective functions 
the values of \newdc and of \newpc agree, and the latter value  is attained, when it is finite. System \eref{dform} can be recast in the primal form 
\beq \label{By0} 
\B x \leq_{K^*} y_0, 
\eeq
where $\B$ and $y_0$ satisfy 
$\R(\B) = \N(\A^*)$ and $\A^*y_0 = c.$ 
It is straightforward 
to show that \eref{dform} is well behaved, if and 
only if \eref{By0} is, and to translate the 
conditions of Theorem \ref{unif-d-thm} 
to characterize when \eref{dform} is well- or badly behaved. 
We leave the details to the reader.

In the special case of semidefinite systems we can obtain the 
following result:  
\bth
Suppose that in the system 
\beq \label{dualsdp}
Y \succeq 0, \, A_i \bullet Y \, = \, c_i \, (i=1, \dots, m)
\eeq
the maximum rank feasible matrix is 
\beq \nonumber
\bar{Y} \, = \, \bpx I_r & 0 \\ 0 & 0 \epx \, {\rm for \,  some} \, r \geq 0. 
\eeq
Then \eref{dualsdp} is badly behaved if and only if there is a matrix $V$ 
and  a real number $\lambda$ 
such that 
\beq \nonumber 
A_i \bullet V \, = \, \lambda c_i \,\, (i=1, \dots, m),
\eeq
and 
$$
V \, = \, 
\bpx V_{11} & V_{12} \\
                V_{12}^T & V_{22} 
\epx, \, 
$$
where $V_{11}$ is $r$ by $r, \,$ $V_{22} \succeq 0, \,$ and 
$\R(V_{12}^T) \not \subseteq \R(V_{22}).$ 
\qed \enth

We can apply similar arguments to conic linear systems in a subspace form
$$
K \cap (L + x_0),
$$
to characterize their well- or badly behaved status.

We can also characterize badly behaved second order conic systems 
similarly as we did it for \eref{p-sd} in Theorem \ref{badsdp}. 
This result is in version 2 of the online version of the paper on arxiv.org.

We finally mention a subject for possible future work.  
The interplay of algebraic geometry and optimization is an 
active research area: see for instance the recent monograph by Blekherman et al  \cite{Blekhetal:12}, and the paper 
of Klep and Schweighofer \cite{KlepSchw:13}. 
It would be interesting to see how our certificates 
of bad and good behavior can be interpreted in the language of algebraic geometry. 

\co{
We finally mention a subject for possible future research:  interpreting our certificates 
of bad and good behavior in the language of algebraic geometry. 
The interplay of algebraic geometry and optimization is an 
active research area: see for instance the recent monograph \cite{Blekhetal:12}. It would be interesting to see how our certificates 
of bad and good behavior can be interpreted in the language of algebraic geometry. }

\appendix
\section{A larger badly behaved semidefinite system}

In this appendix we give  a larger badly behaved semidefinite system to illustrate the standard form reformulation (\ref{p-sd-bad}).  
What is nice about this example is that the bad behavior of the 
original (not reformulated) system is very difficult to 
verify by an {\em ad hoc} argument, whereas the bad behavior of the reformulated system 
is self-evident.

\bex{\rm \label{big-ex}
Consider the badly behaved semidefinite system
\begin{small}
\beq \label{messyex4}
\ba{l} 
x_1  \bpx 4  &   3 &   -5  &  -3 \\
3 &   -2 &    0 &   -2 \\
-5 &    0 &  -12 &   -8 \\
-3 &   -2 &   -8 &   -4
\epx 
+ x_2  \bpx 14   & 10 &  -15 &   -9 \\
10 &   -6 &    0 &   -6 \\
-15 &    0 &  -36  & -24 \\
-9 &   -6  & -24  & -12 \epx 
+ x_3  \bpx   8 &    6 &   -5 &   -3 \\
6 &   -4  &   0 &   -2 \\
-5  &   0  & -12  &  -8 \\
-3 &   -2 &   -8 &   -4
\epx
+ x_4  \bpx 20   & 15 &  -25 &  -13 \\
15 &  -10 &   -1 &   -9 \\
-25 &   -1 &  -58  & -38 \\
-13 &   -9 &  -38 &  -18
 \epx \\ \preceq 
\bpx  45  &  32 &  -55 &  -31 \\
32 &  -19 &   -1 &  -21 \\
-55 &   -1 & -130 &  -86 \\
-31 &  -21 &  -86 &  -42 \epx.
\end{array}
\eeq
\end{small}      

We show how to bring \eref{messyex4} into the form of \eref{p-sd-bad}, so let us 
denote the constraint matrices on the left by $A_i \, (i=1, \dots, 4), \,$ 
and the right hand side matrix by $B. $ Let 
$$
T \, = \, \bpx 1 & 0 & 0 & 0 \\
                    0 & 1 & 0 & 0 \\
                     0 & 0 & -1/2 & 1/2 \\
                     0 & 0 & 3/2 & -1/2 \epx,
$$        
apply the rotation $T^T()T$ to all $A_i$ and $B, \,$ then perform the following 
operations:
\beqast
B & := & B - A_1 - 2A_2 + A_3 - A_4, \\
A_4 & := & -5A_1 + A_4, \\
A_3 & := & -2A_1 + A_3, \\
A_2 & := & -3A_1 + A_2, \\
A_1&  := & A_1 - 2A_2 + A_3. 
\eeqast
We obtain the system 
\begin{small}
\beq \label{cleanex5} 
\ba{l} 
x_1 \bpx 
     0 &   1  &    0 &    0 \\
     1 &   -2 &    0 &    0 \\
     0 &    0 &    0 &    0 \\
     0 &    0 &    0 &    0 \epx
+ x_2 \bpx 
     2 &    1 &    0 &    0 \\
     1 &    0 &    0 &    0 \\
     0 &    0 &    0 &    0 \\
     0 &    0 &    0 &    0 \epx
+ x_3 \bpx
     0  &   0 &    2 &    1 \\
     0  &   0 &    3 &   -1 \\
     2  &   3 &    0 &    2 \\
     1  &  -1 &    2 &    0 \epx
+ x_4 \bpx
     0 &    0 &    3 &   -1 \\
     0 &    0 &    2 &   -1 \\
     3 &    2 &    2 &    0 \\
    -1 &   -1 &    0 &    0 \epx \\
 \preceq 
\bpx 
     1 &     0 &     0 &     0 \\
     0 &     1 &     0 &     0 \\
     0 &     0 &     0 &     0 \\
     0 &     0 &     0 &     0 \epx.
\end{array} 
\eeq
\end{small}
In (\ref{cleanex5}) the matrices 
\beq \label{y1y2}
Y_1 \, = \,  \bpx 0  &   0 &    0  &   0 \\
0  &   0 &    0  &   0 \\
0  &   0 &    0  &   0 \\
0  &   0 &    0  &   1 
\epx, \, {\rm and} \, Y_2 \, = \,  \bpx 0  &   0 &    0  &   1 \\
0  &   0 &    0  &   1 \\
0  &   0 &    2  &   0 \\
1  &   1 &    0  &   0 
\epx
\eeq
are orthogonal to all the constraint matrices, thus they prove that the rank of any 
slack is at most two. So in (\ref{cleanex5}) 
the right hand side is the maximum rank slack. 

It is easy to see that (\ref{cleanex5}) is badly behaved: following the proof of the ``If'' implication in 
Theorem \ref{badsdp-re}, one can see that the objective 
function $ \sup - x_4$ yields a value of $0$ over \eref{cleanex5}, but there is no 
 dual solution with the same value.  
}
\eex

\co{We obtain an interesting example 
by slightly modifying Tun\c{c}el's example on page 43 in the recent book \cite{Tuncel:11}.
\bex \label{levent-ex} 
{\rm For a positive integer $n$ let us define the symmetric unit matrices
$E_i = e_i e_i^T, \,$ and $E_{ij} = e_i e_j^T + e_j e_i^T, \,$ if $i \neq j.$ 

Let $n \geq 4, \,$ and define the system 
\beq \label{levent-sys}
\sum_{i=2}^{n-1} x_i A_i \preceq B,
\eeq
where 
$A_i = E_{i} + E_{1,i+1} \, (i=2, \dots, n-1), \,$ and $B = E_{11}.$ 
(We number the matrices from $2$ to $n-1$ for convenience.)

This system is badly behaved. 
First note that the matrices 
$$ 
Y_1 \, = \, E_n, \, Y_i \, = \, 2 E_{n-i+1} - E_{1, n-i+2} (i=2, \dots, n-1)
$$
prove that $B$ is a maximum rank slack, and each of the $A_i$ matrices can serve as a $V$ matrix 
of Theorem \ref{badsdp}. This system is in the form of \eref{p-sd-bad}. 

If we change the $A_i$ to $A_i = E_{i} - E_{i+1} + E_{1,i+1}, \,$ then 
the system \eref{levent-sys} becomes well behaved. It is in the form 
of \eref{p-sd-good}, with $Y_1 = 0 \oplus I_{n-1}$ proving that $B$ is a maximum rank slack.
}
\eex

}
\section{Proof of Lemmas \ref{duffin} and  \ref{psd-dir}}

In this section we prove Lemmas \ref{duffin} and 
\ref{psd-dir}.

First we need some definitions and notation. 
For optimization problems we use the symbol $\val()$ to denote their optimal value. 
For program $(\dc)$ we say that  
$\{ y_i \} \subseteq K^*$ is an {\em asymptotically feasible (AF)} solution,
if $ \A^* y_i \rightarrow c, $ and the {\em asymptotic value of $(\dc)$} is 
$$
\ba{rcl}
\aval(\dc) & = & \inf \{ \, \lim \, b^* y_i \, | \, \{ y_i \} \, {\rm is \; asymptotically \; feasible \; in  \;\;} (\dc) \, \},
\ena
$$ 
where the infimum is taken over those AF solutions for which 
$\lim \, b^* y_i$ exists. 

We prove Lemma \ref{duffin} by adapting an argument from
\cite{DufJerKar:81}.  We also  rely on the following lemma due to Duffin:
\ble \label{asd} (Duffin \cite{Duff:56}) 
Problem $(\pc)$ is feasible with $\val(\pc) < + \infty, \,$ iff 
$(\dc)$ is asymptotically feasible with $\aval(\dc) > - \infty, \,$ and if these equivalent statements hold,  
then 
\beq \nonumber 
\val(\pc) = \aval(\dc).
\eeq
\qed
\ele

\pf{of Lemma \ref{duffin}} 
We will use the notation 
$$
\A_h = \bpx A & b \\ 0 & 1 \epx
$$
(which is also used in the proof of Theorem \ref{unif-d-thm}). 

\pf{(If)} Suppose that $\A_h^*(K \times \rad{}_+)^*$ is closed and let $c$ be an objective vector, such that $c_0 := \val(\pc)$ is finite. 
Then $\aval(\dc) = c_0 \,$ 
holds by Lemma \ref{asd}, so there is
$
\{ y_i \} \subseteq K^* \; \mathrm{s.t.} \; \A^* y_i  \rightarrow  c, \; \mathrm{and} \; b^* y_i  \rightarrow  c_0, \; 
$
i.e., 
\beq \nonumber 
(c, c_0) \in \, \cl (\A,b)^* K^* \, \subseteq  \, \cl \A_{h}^* (K^* \ti \rad{}_+) \, = \, \A_{h}^* (K^* \ti \rad{}_+).
\eeq
Hence there is  $y \in K^*, \, s \geq 0$ such that $\A^* y = c, \, $ and $b^* y + s = c_0;$ by weak duality
$b^* y = c_0$ must hold. So $y$ is a feasible solution of $(\dc)$  with value $c_0, \,$ and 
this completes the proof. 

\pf{(Only if)}  
To obtain a contradiction, suppose that 
$\A_h^*(K \times \rad{}_+)^*$ is not closed;  then we will show 
that \eref{p} is badly behaved. Let us choose $c$ and $c_0$ such that  
\beq \nonumber 
\begin{array}{rclrcl}
	(c, c_0) & \in & \cl \A_{h}^* (K^* \ti \rad{}_+) \setminus \A_{h}^* (K^* \ti \rad{}_+).
\end{array}
\eeq
By $(c, c_0) \in \cl \A_{h}^* (K^* \ti \rad{}_+)$ there is 
$\{ (y_i, s_i) \} \subseteq K^* \ti \rad{}_+ \; \mathrm{s.t.} \;
\A^* y_i  \rightarrow  c, \; \mathrm{and} \; b^* y_i + s_i  \rightarrow  c_0. $ 
Hence
$$\val(\pc) \, = \, \aval(\dc)   \leq  c_0, \,$$ 
where the equality comes from Lemma  \ref{asd}.  

However, $(c, c_0) \not \in \A_{h}^* (K^* \ti \rad{}_+)$ shows that no feasible solution of $(\dc)$ can have 
value $\leq c_0$. Hence either $\val(\dc)  >  c_0 \,$ (this includes the case $\val(\dc) = + \infty, \,$ 
i.e., when $(\dc)$ is infeasible), or $\val(\dc)$ is not attained.
\qed

To prove Lemma \ref{psd-dir} we need another lemma, which is mostly based on results 
surveyed in \cite{Pataki:00A}. 
\ble \label{ldir}
Let $C$ be a closed convex cone, $x \in C, \,$ and $E$ the smallest face of $C$ that contains 
$x.$ Then 
\beqa \label{lemma1-dir}
\dir(x, C) & = & C + \lin E,  \\ \label{lemma1-ldir}
\ldir(x, C) & = & \lin E, \\ \label{lemma1-cldir}
\cl \dir(x, C) & = & (C^* \cap x^\perp)^*, \\ \label{lemma1-tan}
\tan(x, C)  & = & (C^* \cap x^\perp)^\perp.
\eeqa
\ele
\pf{} Statements \eref{lemma1-dir} and \eref{lemma1-cldir} 
are in Lemma 3.2.1 in \cite{Pataki:00A} (in Lemma 2.7 in the online version).
We also proved statement \eref{lemma1-tan} 
there, assuming that $C$ is nice. In fact, it 
follows from \eref{lemma1-cldir} and \eref{deftan} in general. 

In \eref{lemma1-ldir} the containment $\supseteq$ is trivial. To see 
$\subseteq$ let $y \in \ldir(x,C), \,$ then 
$x \pm \eps y \in C$ for some $\eps > 0.$ 
Hence $x \pm \eps y \in E, \,$ so $\eps y \in \lin E, \,$ and this completes the 
proof. 
\qed

\pf{of Lemma \ref{psd-dir}}
Let $F$ be the smallest face of $\psd{n}$ that contains $Z.$ 
Then clearly $F = \psd{r} \oplus \{  0  \}, \,$ and 
$\psd{n} \cap Z^\perp = \{  0 \, \} \oplus \psd{n-r}.$ 
Hence statements \eref{psdldir}-\eref{psdtan} follow by taking 
$C = \psd{n}, \, x = Z, \, E = F \,$ in Lemma \ref{ldir}. 

Next, fix $Y \in \cl \dir(Z, \psd{n}), \,$ and partition it as in the right hand side set in 
\eref{psdcldir}.
Then \eref{psdcldirdir} is equivalent to 
\beq \label{toprovedir}
Y \in \dir(Z, \psd{n}) \, \LRA \, \R(Y_{12}^T) \subseteq \R(Y_{22}).
\eeq
Let $P$ be an orthogonal matrix, such that 
$P^T Y_{22}P = I_s \oplus 0, \,$ where 
$s$ is the number of positive eigenvalues of $Y_{22}$ and $T = I_r \oplus P.$ 

Define
$$
V \, := \, T^T Y T \, = \, \bpx Y_{11} & Y_{12} P \\
               P^T Y_{12}^T & P^T Y_{22} P 
                 \epx \, = \, \bpx Y_{11} & Y_{12} P \\
               P^T Y_{12}^T & I_s \oplus 0  
                 \epx.
$$
Next we claim 
\beqa \label{eq1}
Y \in \dir(Z, \psd{n}) & \LRA & V \in \dir(Z, \psd{n}), \,\, \\ \label{eq2} 
\R(Y_{12}^T) \subseteq \R(Y_{22}) & \LRA & \R(P^T Y_{12}^T) \subseteq \R(P^T Y_{22} P).
\eeqa
Indeed, \eref{eq1} follows 
from $T^T Z T = Z, \,$ and the definition of feasible directions.
As to \eref{eq2}, the left hand side statement 
holds, iff there is a matrix $D$ with
\beq \label{eq11}
Y_{12}^T \, = \, Y_{22} D,
\eeq
and the right hand side statement holds, iff there is a matrix $D'$ such that 
\beq \label{eq21}
P^T Y_{12}^T \, = \, P^T Y_{22} P D'.
\eeq
If $D$ satisfies \eref{eq11}, then 
$D' := P^{-1}D$ satisfies \eref{eq21}.
Conversely, if \eref{eq21} holds for $D',$ then 
$D := P D'$ verifies \eref{eq11}. 

Next, partition $Y_{12}P$ as $(V_{12}, V_{13}), \,$ 
so that  $V_{12}$ has $s$ columns; then \eref{eq2}  
is equivalent to $V_{13} = 0.$ So we
only need to prove
\beq \label{toprovedir-R}
V \in \dir(Z, \psd{n}) \, \LRA \, V_{13} = 0.
\eeq
Consider the matrix $Z + \eps V \,$ for some $\eps > 0.$ 
If $V_{13} \neq 0, \,$ 
then $Z + \eps V$ is not positive semidefinite for any $\eps>0, \,$ and 
this proves the direction $\RA.$ As to  $\LA, \,$ 
if $V_{13} = 0, \,$ then 
by the Schur-complement condition for positive semidefiniteness we have 
that $Z + \eps V \succeq 0$ iff 
$$
(I_r + \eps V_{11}) - (\eps V_{12})(\eps I_s)^{-1}(\eps V_{12}^T) \succeq 0,
$$
and the latter is clearly true for some small $\eps > 0.$ 
\qed

{\bf Acknowledgement} My sincere thanks are due to the anonymous referees, the Associate Editor, and Shu Lu for their careful reading of the manuscript, 
and thoughtful  comments. I also thank Minghui Liu for helpful comments, and for his 
help in proving Theorem \ref{goodsdp-re}. My thanks are also due to Asen Dontchev 
for his support while writing  this paper.


\bibliography{mysdp}

\end{document}